\newfont{\Bbb}{msbm10 scaled\magstephalf}
\documentclass[12pt]{amsart}
\usepackage{amsmath, amsthm, amscd, amsfonts}
\newcommand{\WW}{{\mathbb W}}
\newcommand{\RR}{{\mathbb R}}

\newcommand{\CC}{{\mathbb C}}

\newcommand{\NN}{{\mathbb N}}

\newcommand{{\Li}}{{\mathcal L}_p}
\newcommand{{\Bl}}{{\mathcal B}_p}

\newtheorem{theorem}{Theorem}[section]
\newtheorem{lemma}[theorem]{Lemma}

\newtheorem{corollary}[theorem]{Corollary}

\newtheorem{proposition}[theorem]{Proposition}

\def\msk{\medskip}
\def\ol{\overline}
\def\bege{\begin{equation}} \def\ende{\end{equation}}
   \def\om{\omega}
 \def\g{\gamma}\def\ol{\overline} \def\ve{\varepsilon} \def\d{\delta}
\def\pt{\partial} 
\def\begr{\begin{eqnarray}} \def\endr{\end{eqnarray}}
\def\qand{\quad\mbox{ and }\quad}

\def\bnum{\begin{enumerate}}  \def\enum{\end{enumerate}}

\begin{document}

\title[Composition operators on generalized Bloch spaces]{COMPOSITION OPERATORS BETWEEN
GENERALIZED BLOCH SPACES OF THE POLYDISK}

\author[D. Clahane, S. Stevi\'c, and Z. Zhou]{Dana D. Clahane$^*$ \and Stevo Stevi\'c \and Zehua Zhou$^+$}
\address{\newline Department of Mathematics\newline Indiana
University\newline
Bloomington, IN 47405.
\newline{\em Current
Address:}
\newline Department of Mathematics\newline
University of California
\newline Riverside, CA 92521.}

\email{dclahane@math.ucr.edu}
\address{\newline Mathematical Institute of Serbian Academy of Science\newline Knez
Mihailova 35/I\newline 11000 Beograd, Serbia.}
\email{sstevic@ptt.yu; sstevo@matf.bg.ac.yu}

\address{\newline Department of Mathematics
\newline Tianjin University\newline Tianjin 300072, China.}
\email{zehuazhou@hotmail.com; zehuazhou2003@yahoo.com.cn}



\date{}
\thanks{\noindent $^*$Work by the first author was supported by
an NSF VIGRE Postdoctoral Fellowship at Indiana University
Bloomington.\\
\indent $^+$Zehua Zhou, corresponding author. Supported in part by
the National Natural Science Foundation of China (Grant $\#$'s
10001030 and 10371091), and LiuHui Center for Applied Mathematics,
Nankai University and Tianjin University.}

\begin{abstract} Let $p$, $q>0$, and suppose that $\phi$ is a holomorphic
self-map of the open unit polydisk $U^n$ in $\CC^n$.  We extend to
$U^n$ previous one-variable characterization results of K. Madigan
on the $p$-Lipschitz space ${\mathcal L}_p(U^n)$ and K. Madigan/A.
Matheson on the Bloch space ${\mathcal B}^1(U)$ by obtaining
function-theoretic conditions on $\phi$ such that the composition
operator $C_{\phi}$ is bounded or compact between the generalized
Bloch spaces ${\mathcal B}^p(U^n)$ and ${\mathcal B}^q(U^n)$.  These
conditions turn out to be different in the cases $p\in(0,1)$ and
$p\geq 1$.  We also obtain corresponding characterization results
for composition operators between generalized little Bloch spaces
${\mathcal B}_0^p(U^n)$ and ${\mathcal B}_0^q(U^n).$
\end{abstract}

\maketitle

\section{Introduction}

Suppose that $D$ is a domain in $\CC^n$, and let $\phi$ be a
holomorphic map from $D$ to itself.  Denoting the vector space of
holomorphic functions on $D$ by ${\mathcal H}(D)$, we define the
{\em composition operator} $C_\phi$ on ${\mathcal H}(D)$ by $
C_\phi(f)=f\circ\phi$. This operator is linear.

A currently pursued theme in the theory that surrounds
function-theoretically defined operators such as composition
operators is the general problem of characterizing, in terms of
mapping properties of the operator-inducing function, those
operators that have given properties, such as boundedness and
compactness (cf. \cite{Russo composition}).
 In this paper, we further this theme by obtaining function-theoretic
conditions on a holomorphic self-map $\phi$ of the open unit
polydisk $U^n$ in $\CC^n$ such that $C_\phi$ is bounded or compact
between $p$- and $q$-Bloch spaces of $U^n$ in $\CC^n$ for $p>0$ or
$p$- and $q$-Lipschitz spaces of $U^n$ for $0<p<1$. These results
generalize previously known one-variable results obtained by K.
Madigan/A. Matheson \cite{mm2} and K. Madigan \cite{mm1},
respectively.

In questions about composition operators, $C_\phi$ is usually
restricted to a subspace of ${\mathcal H}(D)$, such as a weighted
Hardy, weighted Bergman, Lipschitz, or Bloch space over a given
domain $D$ in $\CC^n$.  These restricted composition operators often
have properties that are encoded purely and beautifully in the
behavior of $\phi$ on $D$.

Despite the simplicity of the definition of $C_\phi$, it is not
uncommon for solutions of problems involving this type of operator
to require profound and interesting analytical machinery; moreover,
the study of composition operators has arguably become a major
driving force in the development of modern complex analysis.  The
lucent texts \cite{Cowen/MacCluer}, \cite{Shapiro}, \cite{zh}, and
conference proceedings collection \cite{cmr} are good sources for
information about many of the developments in the theory of
composition operators up to the middle of the last decade.

The theory of composition operators on classical function spaces
such as Hardy and Bergman spaces of the unit disk $U$ is now highly
advanced, but the multivariable situation remains mysterious.  The
Littlewood Subordination Principle
(\cite[~Thm.~2.22]{Cowen/MacCluer}, \cite{Littlewood}) leads to the
foundational result that the restriction of $C_\phi$ to the Hardy or
standard weighted Bergman spaces of $U$ is always bounded. J.
Cima/W. Wogen confirmed in \cite{Cowen/MacCluer} that the
multivariable situation is much different by making the pivotal
observation that composition operators can be unbounded on the Hardy
space of the open unit ball in $\CC^2$.  These operators can also be
unbounded on many other functional Banach spaces over domains in
$\CC^n$, even when $n=1$ (cf. \cite{mm2}). Therefore, in such a
setting it is natural to pursue the objective of
function-theoretically or geometrically characterizing those maps
$\phi$ that induce bounded or compact composition operators.

An intriguing aspect of the corresponding results for composition
operators on spaces such as these, which are defined by derivative
rather than integral growth, is that many one-variable results about
composition operators on such spaces generalize to the several
variable case. In addition, Lipschitz and Bloch spaces share duality
relationships with the classical function spaces, so understanding
composition operators on Lipschitz and Bloch spaces could shed light
on the corresponding open problems for Hardy and weighted Bergman
spaces.

In \cite{mm1}, K. Madigan determined the analytic self-maps $\phi$
of $U$ for which $C_\phi$ is bounded on the analytic $p$-Lipschitz
spaces and $(1-p)$-Bloch spaces of $U$ for $0<p<1$. His
characterization, which we state below, uses a classical
Hardy/Littlewood result \cite{hl} that essentially asserts equality
of the $p$-Lipschitz spaces ${\Li}(U)$ and $(1-p)$-Bloch spaces of
$U$, along with the stronger fact that the natural norms on these
spaces are equivalent (see Lemma \ref{hln} in the present paper):
\begin{theorem}\label{mot}
Suppose that $p\in(0,1)$, and let $\phi$ be a holomorphic self-map
of $U$. Then $C_\phi$ is bounded on ${\Li}(U)$ if and only if
\[
\sup_{z\in
U}\left(\frac{1-|z|^2}{1-|\phi(z)|^2}\right)^{1-p}|\phi'(z)|<\infty.
\]
\end{theorem}
Shapiro \cite{Shapiro small} proved that on
automorphism-invariant, boundary-regular small spaces of $U$,
compactness of $C_\phi$ implies that $||\phi||_\infty<1$. Such
spaces include ${\Li}(U)$ for $0<p\leq 1$, and Shapiro's result
holds more generally on the unit ball $B_n=\{ z\, |\, |z|<1\}$ as
well \cite[~p.~180]{Cowen/MacCluer}.  Known proofs of this fact
depend on the strict convexity of $B_n$, so as one would expect,
we treat the issue of compactness of $C_\phi$ on ${\Li}(U^n)$ with
a different approach.

By Lemma \ref{hln} in the present paper, $\Li(U)$ in Theorem
\ref{mot} can be replaced by ${\mathcal B}^{1-p}(U)$. It is
therefore natural to wonder what results can be proven about
boundedness and compactness of $C_\phi$ on $p$-Bloch spaces for
arbitrary positive numbers $p$ or, more generally, between possibly
different $p$- and $q$-Bloch spaces of multivariable domains.  In
this paper, we answer these questions completely for $U^n$, although
we have opted in this paper not to treat boundedness and compactness
with a unified essential norm approach. In fact, the corresponding
results for essential norms of composition operators can be
routinely formulated using the results of this paper as a template.

In \cite{mm2}, Madigan/A. Matheson showed that $C_\phi$ is bounded
on the Bloch space ${\mathcal B}^1(U)$.  They additionally obtained
the following result:
\begin{theorem}[Madigan/Matheson]\label{madmath}
If $\phi$ is a holomorphic mapping of the unit disk $U$ into
itself, then $\phi$ induces a compact composition operator on
${\mathcal B}^1(U)$ if and only if for every $\ve>0$ there is an
$r\in(0,1)$ such that whenever $|\phi(z)|>r$,
$$\frac{1-|z|^2}{1-|\phi(z)|^2}|\phi'(z)|<\ve.$$
\end{theorem}

 The analogues of these facts for bounded symmetric domains were
obtained by Z. Zhou/J. Shi in \cite{zs3}.  In this paper, we
characterize bounded and compact composition operators between
$p$-Bloch and $q$-Bloch spaces of $U^n$ for all values of $p$ and
$q$, thus extending Theorems \ref{mot} and \ref{madmath} in several
ways. We also obtain the corresponding results for composition
operators between suitably defined ``little" $p$- and $q$-Bloch
subspaces ${\mathcal B}^p_0(U^n)$ and ${\mathcal B}^q_0(U^n)$.

\section{Notation and Additional Motivating Results}

The reader is referred to \cite{Rudin} for standard several complex
variables terminology.  Letting $n\in\NN$, we denote the {\em open
unit polydisk} in $\CC^n$ by $U^n$.  $U^n$ is the domain consisting
of all points $(z_1,z_2,\ldots,z_n)$ whose coordinates $z_i$ for
$i=1,2\ldots,n$ are in the open unit disk $U$ in $\CC$.

Denote the set of positive integers by $\WW$.  $\gamma$ is said to
be an $n$-{\em multi-index} if and only if $\gamma\in \WW^n$, where
$\WW$ denotes the set of whole numbers. If
$\gamma=(\gamma_1,\gamma_2,\ldots,\gamma_n)$ is an $n$-multi-index,
then we define $z^\gamma$ for each $z\in\CC^n$ to be
$z_1^{\gamma_1}z_2^{\gamma_2}\ldots z_n^{\gamma_n}$, with the
convention that $z_i^0=1$ for all $i\in\{1,2,\ldots,n\}$ in this
context.

For $z,w\in\CC^n$, define $[z,w]_j=z$ if $j=0$, $[z,w]_j=w$ if
$j=n$, and
\[
[z,w]_j=(z_1,z_2,\ldots,z_{n-j},w_{n-j+1},\ldots w_n)
\]
if $j\in\{1,2,\ldots,n-1\}$.  Then $[z,w]_{n-j+1}=w$ when $j=1$,
and $[z,w]_{n-j+1}=z$ when $j=n+1$.  For $j=2,3\ldots n$, note
that
\[
[z,w]_{n-j+1}=(z_1,z_2,\ldots,z_{j-1},w_j,\ldots w_n).
\]
If $a\in\CC$, then we write
\[z_j'=(z_1,z_2,\ldots ,z_{j-1},z_{j+1},\ldots ,z_n)\]
and
\[(a,z_j')=(z_1,z_2,\ldots,z_{j-1},a,z_j,\ldots,z_n).
\]
If $A\subset\CC^n$, $z\in\CC^n$, and $d(z,A)$ is the Euclidean
distance from $z$ to $A$, then we abbreviate the condition that
$d(z,A)\rightarrow 0$ by $z\rightarrow A$.

Let $(z^{[j]})_{j\in\NN}$ be a sequence in $\CC^n$.  We will adopt
the notation $(z^{[j']})_{j\in\NN}$ for an arbitrary subsequence of
$(z^{[j]})_{j\in\NN}$; that is, the notation means that we have a
self-map $j\rightarrow j'$ of $\NN$ such that $j'<(j+1)'$ for all
$j\in\NN$.

For $p\in(0,1]$, the holomorphic {\em $p$-Lipschitz space}
${\Li}(U^n)$ consists of all $f\in {\mathcal H}(U^n)$ such that for
some $M\in [0,\infty)$ and all $z,w\in U^n$,
\[
|f(z)-f(w)|\leq M|z-w|^p.
\]
The norm that we will naturally associate with ${\Li}(U^n)$ is given
by
\[
||f||_{{\Li}}:=|f(0)|+\sup_{z\not=w\in
U^n}\frac{|f(z)-f(w)|}{|z-w|^p}.
\]
${\Li}(U^n)$ is a Banach space with respect to this norm.

For an example of a holomorphic map $\phi:U\rightarrow U$ such that
$C_\phi$ is not bounded on ${\Li}(U)$, see \cite{mm1}.
 $\psi:U^n\rightarrow U^n$ defined by $\psi(z_1,z_2,\ldots
,z_n)=(\phi(z_1),0,0,\ldots ,0)$ then induces an unbounded
composition operator $C_\psi$ on ${\Li}(U^n)$ and ${\Li}(B_n)$.

A remark about the case $\alpha=1$ is in order here.  As is the case
in one variable \cite[~p.~186]{Cowen/MacCluer}, characterizing the
bounded composition operators on ${\mathcal L}_1(U^n)$ is not
difficult, so we leave the proof of the following statement as an
exercise for the interested reader:

\medskip

\indent{\em Let $\phi:U^n\rightarrow U^n$ be holomorphic.  Then
$C_\phi$ is bounded on ${\mathcal L}_1(U^n)$ if and only if
$\phi_j\in {\mathcal L}_1(U^n)$ for all $j=1,2,\ldots n$.}

\medskip

The so-called {\em Bloch space} ${\mathcal B}(U^n)$ can be defined
as follows.  The following formula for the Bergman metric
$H:U^n\times\CC^n\rightarrow [0,\infty)$ on $U^n$ appears, for
example, in \cite{t}:
\begin{equation}\label{berg}
H(z,u)=\sum_{k=1}^n\frac{|u_k|^2}{(1-|z_k|^2)^2}.
\end{equation}
We can then characterize the Bloch space ${\mathcal B}(U^n)$ as the
vector space of all $f\in {\mathcal H}(U^n)$ such that
\[
b_1(f)=\sup_{z\in U^n}Q_f(z)<\infty,
\]
where

$$Q_f(z)=\sup_{u\in \CC^n\setminus\{0\}}\;\frac{|\langle\nabla f(z),
u\rangle|}{\sqrt{H(z,u)}}\,{\text{ and }}\,\nabla
f(z)=\left(\frac{\pt f}{\pt z_1}(z),\ldots,\frac{\pt f}{\pt
z_n}(z)\right).$$

One can replace $U^n$ in the definition of $b_1(f)$ by more
general homogeneous domains, thus leading to a definition of Bloch
space on these domains.  As is the case on $U$, Bloch functions on
homogeneous domains have a variety of characterizations that
connect Bloch spaces to several cursorily unexpected concepts in
analysis \cite{t}.

Using \cite[Inequ.~(1.6),~p.~3]{t2}, one can verify that the
quantities $$|f(0)|+b_1(f)\,\,\,\,\,{\text{ and }}\,\,\,\,\,
\|f\|_1=|f(0)|+\sup\limits_{z\in
U^n}\sum\limits^n_{k=1}\left|\frac{\partial f}{\partial
z_k}(z)\right|(1- |z_k|^2)$$ define equivalent norms on ${\mathcal
B}(U^n)$.  This is a convenient starting point for our introduction
of a generalized Bloch space, which we will call a {\em $p$-Bloch
space}.

Let $p>0$.  We say that $f\in {\mathcal H}(U^n)$ belongs to the
$p$-Bloch space ${\mathcal B}^p(U^n)$ if and only if for some
$M\in[0,\infty)$ and all $z\in U^n$, we have that
\[
\sum\limits^n_{k=1} \left|\frac{\partial f} {\partial
z_k}(z)\right|\left(1- |z_k|^2\right)^p\leq M.
\]
A natural norm that makes ${\mathcal B}^p(U^n)$ a Banach space is
given by
\[
\|f\|_{\Bl}=|f(0)|+\sup\limits_{z\in U^n}\sum\limits^n_{k=1}
\left|\frac{\partial f} {\partial z_k}(z)\right|\left(1-
|z_k|^2\right)^p.
\]
${\mathcal B}^p(U^n)$ with the above-given norm is a Banach space.

In order to define a meaningful ``little" version of these
$p$-Bloch spaces, we should first point out a significant
difference between Bloch spaces of $U$ and $U^n$ for $n>1$.
${\mathcal B}^1_0(U)$ is customarily defined (cf.
\cite[~p.~83]{zh}) to be the closed subspace of functions $f\in
{\mathcal B}^1(U)$ such that
$$\lim_{z\to\pt
U}|f'(z)|(1-|z|^2)=0.$$  ${\mathcal B}^1_0(U)$ turns out to be the
closure of the polynomials in ${\mathcal B}^1(U)$.  The natural
extensions of these two characterizations of Bloch spaces do not
coincide on $U^n$ for $n>1$ (nor for any bounded symmetric domain
of rank greater than $1$ \cite{t2}); in fact, for the reader's
convenience, we modify ideas given in \cite[~p.~14]{t2} to prove
the following statement:
\begin{proposition}\label{constant}
Let $f\in {\mathcal H}(U^n)$, and assume that $p>0$.  Then
$$\lim_{z\to\pt U^n}\sum\limits^n_{k=1}\left|\frac{\partial f}{\partial
z_k}(z)\right|(1- |z_k|^2)^p=0,$$ if and only if $f$ is constant.
\end{proposition}
\begin{proof}
Fixing $a\in U$, one sees that $g:U^{n-1}\rightarrow\CC$
defined by
\[g(z')=\frac{\pt f}{\pt z_1}(a,z')(1-|a|^2)^p\]
is holomorphic, where $z'=(z_2,z_3,\ldots,z_n)\in U^{n-1}.$  If
$z'\to \pt U^{n-1}$, then certainly, $(a,z')\to \pt U^{n},$ which
implies by hypothesis that
$$\lim_{z'\to \pt U^{n-1}}\left|\frac{\pt f}{\pt z_1}(a,z')\right|(1-|a|^2)^p=0.$$
By the maximum principle, we have that for every $z'\in U^{n-1}$,
$$\frac{\pt f}{\pt z_1}(a,z')(1-|a|^2)^p= 0.$$
Since $a\in U$ is arbitrary, it follows that ${\pt f}/{\pt z_1}=0$
on $U^n.$ Repeating this argument by placing $a\in U$ in each
component of $U$, we obtain that ${\pt f}/{\pt z_j}=0$ on $U^n$ for
each $j\in \{2,\ldots,n\}.$  Hence, $f$ is constant. The proof of
the converse is immediate.
\end{proof}
Since Proposition \ref{constant} implies that replacement of the
supremum in the definition of ${\mathcal B}^p(U^n)$ by a little-oh
condition as $z$ tends to $\partial U^n$ leads to an uninteresting
space, we will define the {\em little $p$-Bloch space ${\mathcal
B}_0^p(U^n)$} to be the closure of the polynomials in ${\mathcal
B}^p(U^n)$.

Prior to obtaining the generalization of Theorem \ref{madmath} to
classical bounded symmetric domains and the boundedness of $C_\phi$
on the Bloch spaces of those domains in \cite{zs3}, Z. Zhou/J. Shi
had already shown in \cite{zs1} that $C_\phi$ is always bounded on
${\mathcal B}(U^n)$ and had furthermore formulated a correct
characterization of the compact composition operators on ${\mathcal
B}(U^n)$ as follows:

\begin{theorem}\label{unproven}
Let $\phi$ be a holomorphic self-map of $U^n.$  Then $C_{\phi}$ is
compact on ${\mathcal B}(U^n)$ if and only if for every $\ve>0$
there exists a $\d\in(0,1)$ such that for all $z\in U^n$ with
$\;dist(\phi(z),\pt U^n)<\d$,
$$\sum\limits^n_{k,l=1}\left|\frac{\partial \phi_{l}}
{\partial z_k}(z)\right|\frac{1-|z_k|^2} {1-|\phi_l(z)|^2}<\ve.$$

\end{theorem}
Though the proof in \cite{zs1} of sufficiency is complete, the proof
of necessity contains at least one significant gap.  More
specifically, in view of the reduction assumption $(9)$ in that
proof, one cannot omit consideration of the case when
$|\phi_1(z^{(j)})|\not\to 1$ as $j\to\infty$ for the subsequence
$(z^{(j)})_{j\in U^n} $ that is said there to satisfy
$\phi(z^{(j)})\to \pt U^n$. Our results consider not only
composition operators on ${\mathcal B}^1(U^n)$ but more generally,
composition operators between ${\mathcal B}^p(U^n)$ and ${\mathcal
B}^q(U^n)$; moreover, our proofs omit no cases.

In \cite{z4}, Zhou stated and proved the corresponding compactness
characterization for ${\mathcal B}^p(U^n)$ for $0<p<1$ in an attempt
to see what was needed to overcome the gap in the proof of necessity
in \cite{zs1}; however, the test functions used in \cite{z4} are
only suitable for handling the case $0<p<1$.  As we will see, an
extension of Zhou's characterization in \cite{z4} from $0<p<1$ to
$p=1$ and beyond requires a different family of test functions,
which are defined as follows:

Suppose that $p>0$, $l\in\{1,2,\ldots,n\}$, and $w\in U$.  Then we
define $f^{(l)}_w:{\overline U}^n\rightarrow\CC$,
$g^{(l)}_w:{\overline U}^n\rightarrow\CC$, and $h^{(l)}_w:{\overline
U}^n\rightarrow\CC$ respectively by
\begin{eqnarray*}
f^{(l)}_w(z)&=&\int_0^{z_l}\frac{dt}{(1-\bar wt)^p},\\
g^{(l)}_w(z)&=&\frac{1-|w|^2}{(1-z_l{\bar w})^p},\,{\text{
and }}\\
h^{(l)}_w(z)&=&(z_1+2)(1-|w|^2)^{p-1}g^{(l)}_w(z)\,\,\,{\text{ for
}}l\not=1.
\end{eqnarray*}
We let ${\mathcal T}^p={\mathcal T}^p_1\cup {\mathcal T}^p_2\cup
{\mathcal T}^p_3$, where
\begin{eqnarray*}
{\mathcal T}_1^p&=&\{f^{(l)}_w:w\in U{\text{ and
}}l\in\{1,2,\ldots, n\}\},\\
{\mathcal T}_2^p&=&\{g^{(l)}_w:w\in U{\text{ and }}l\in\{1,2,\ldots,
n\}\},\,{\text{
and }}\\
{\mathcal T}_3^p&=&\{h^{(l)}_w:w\in U{\text{ and
}}l\in\{2,\ldots,n\}\}.
\end{eqnarray*}

To conclude our background remarks, we mention that related results
about composition operators on Bloch or Lipschitz spaces of $U^n$
and classical bounded symmetric domains in $\CC^n$ appear in
\cite{zs2} and \cite{zs3}, for example.  Results of this type for
$B_n$ are given in \cite[~Ch.~4]{Clahane} and \cite{sh}.

The remainder of the present paper is assembled as follows.  In
Section \ref{results}, we state our main results.  In Section
\ref{lemmas}, we state and prove facts that will be needed in the
proofs of those results, which are proven in Section \ref{proofs}.
Corresponding results for composition operators between little
$p$-Bloch and little $q$-Bloch spaces are stated and proved in
Section \ref{little spaces}.

\section{Main Results}\label{results}

The principal objective of this paper is to prove all of the
following theorems and their corollaries:

\begin{theorem}\label{firsttheorem}
Let $\phi$ be a holomorphic
self-map of $U^n$, and suppose that $p,q>0.$  Then the following
statements are equivalent:
\bnum
\item[{\bf (\ref{firsttheorem}a)}] $C_{\phi}$ is a bounded operator from ${\mathcal B}^p(U^n)$
to
 ${\mathcal B}^q(U^n) $;

\item[{\bf (\ref{firsttheorem}b)}] $C_{\phi}$ is a bounded operator from ${\mathcal
B}^p_0(U^n)$ to ${\mathcal B}^q(U^n) $;

\item[{\bf (\ref{firsttheorem}c)}] There is an $M\geq 0$ such that for all $z\in U^n$,
\[
\sum\limits^n_{k,l=1}\left|\frac{\partial \phi_{l}} {\partial
z_k}(z)\right|\frac{(1-|z_k|^2)^q} {(1-|\phi_l(z)|^2)^p}\leq M.
\]
\enum
\end{theorem}

Combining Theorem \ref{firsttheorem} and Lemma \ref{hln1}, the
latter of which will be stated and proved in Section 4, we obtain
the following generalization of Theorem \ref{mot}:

\begin{corollary}\label{lipresult}
Let $\phi$ be a holomorphic self-map of $U^n$, and suppose that
$0<p,q<1$.  Then the following statements are equivalent: \bnum
\item[{\bf (a)}] $C_{\phi}$ is a bounded operator from ${\mathcal L}_{p}(U^n)$
to
 ${\mathcal L}_{q}(U^n) $;

\item[{\bf (b)}] There is a $C\geq 0$ such that for all $z\in
U^n$,
\[
\sum\limits^n_{k,l=1}\left|\frac{\partial \phi_{l}} {\partial
z_k}(z)\right|\frac{(1-|z_k|^2)^{1-q}} {(1-|\phi_l(z)|^2)^{1-p}}\leq
C.
\]
\enum
\end{corollary}
Theorem \ref{firsttheorem} and Corollary \ref{lipresult} put to
rest the issue of boundedness, and the remaining results in this
paper deal with compactness.
\begin{theorem}\label{secondtheorem}
Let $\phi$ be a holomorphic self-map of $U^n,$ and suppose that $p,
q>0.$ \bnum \item[{\bf (a)}] If Condition {\bf
(\ref{firsttheorem}c)} and
\begin{equation}\label{firstcom}
\lim_{\phi(z)\rightarrow\partial
U^n}\sum\limits^n_{k,l=1}\left|\frac{\partial \phi_{l}} {\partial
z_k}(z)\right|\frac{(1-|z_k|^2)^q}{(1-|\phi_l(z)|^2)^p}
=0.
\end{equation}
hold, then $C_{\phi}$ is a compact operator from ${\mathcal
B}^p(U^n)$ and ${\mathcal B}^p_0(U^n)$ to ${\mathcal B}^q(U^n)$.
\item[{\bf (b)}] Let $p\geq 1$ and $q>0.$ If $C_{\phi}$ is a
compact operator from ${\mathcal B}^p(U^n)$ or ${\mathcal
B}^p_0(U^n)$ to ${\mathcal B}^q(U^n)$, then Condition
$(\ref{firstcom})$ holds. \enum
\end{theorem}
The following application of Theorem \ref{secondtheorem} relates
compactness of $C_\phi$ between $p$-Bloch spaces and behavior of
$\phi$ in the Bergman metric:
\begin{corollary}\label{Jac}
Let $\phi:U^n\rightarrow U^n$ be holomorphic, and assume that
$p\geq 1$ and $q\in(0,1]$.  Suppose also that the following
statement holds:

{\bf (\ref{Jac}a)} There exists a $C>0$ such that for all $z\in U^n$
and $u\in\CC^n$,
$$H_z(u,\bar u)\leq CH_{\phi(z)}\big( J_\phi(z)u,\ol{J_\phi(z)u}\big),$$
where
$$J_\phi(z)=\left[\frac{\pt \phi_j}{\pt z_k}(z)\right]_{1\leq j,k\leq n}$$
denotes the Jacobian matrix of $\phi$ at $z$.  Then

{\bf (\ref{Jac}b)} $C_{\phi}$ is not compact from ${\mathcal
B}^p(U^n)$ nor ${\mathcal B}^p_0(U^n)$ to ${\mathcal
B}^q(U^n)$.\end{corollary} {\em Remarks:

(i)} If $\phi\in Aut(U^n),$ $z\in U^n$, and $u\in \CC^n$, then the
following inequality holds \cite{t}:
$$H_z(u,\bar
u)=H_{\phi(z)}\big(J_{\phi}(z)u,\ol{J_\phi(z)u}\big).$$ Hence by
Corollary \ref{Jac}, we observe that if $\phi \in Aut(U^n)$,
$p\geq 1$, and $q\in(0,1]$, then $C_\phi$ is not a compact
operator from ${\mathcal B}^p(U^n)$ or ${\mathcal B}^p_0(U^n)$ to
${\mathcal B}^q(U^n)$.

{\em (ii)} If $\phi$ is a proper holomorphic self-map of $U^n$,
then a manipulation of the definition of limit leads to the
observation that such a map $\phi$ satisfies Condition ({\bf
\ref{firsttheorem}c}) whenever $\phi$ satisfies Equation
(\ref{firstcom}).  Therefore, one can remove from part {\bf (b)}
of Theorem \ref{secondtheorem} the assumption that Condition ({\bf
\ref{firsttheorem}c}) holds if one assumes that $\phi$ is proper.
Nevertheless, Theorem \ref{secondtheorem} implies the following
function-theoretic characterization of compact composition
operators from ${\mathcal B}^p(U^n)$ or ${\mathcal B}^p_0(U^n)$ to
${\mathcal B}^q(U^n)$ for $p\geq 1$ and $q>0$:
\begin{corollary}\label{firstiff}
Let $\phi$ be a holomorphic self-map of $U^n,$ and suppose that
$p\geq 1$ and $q>0.$ Then the following statements are equivalent:
\bnum \item[{\bf (\ref{firstiff}a)}] $C_{\phi}$ is a compact
operator from ${\mathcal B}^p(U^n)$ or ${\mathcal B}^p_0(U^n)$ to
${\mathcal B}^q(U^n)$; \item[{\bf (\ref{firstiff}b)}] Conditions
${\bf (\ref{firsttheorem}c)}$ and $(\ref{firstcom})$ both hold;
\item[{\bf (\ref{firstiff}c)}] $\phi_l(z)\in {\mathcal B}^q(U^n)$
for every $l\in \{1,2,\ldots,n\}$, and  Condition
$(\ref{firstcom})$ holds. \enum
\end{corollary}
Notice that in Theorem \ref{secondtheorem} and Corollary
\ref{firstiff}, the limiting condition involves approach of
$\phi(z)$ to the boundary of $U^n$ rather than toward the
distinguished boundary, as is the case in Theorem \ref{p unit
interval} below for the case $p\in(0,1)$:
\begin{theorem}\label{p unit interval}
Let $\phi$ be a holomorphic self-map of $U^n$, and assume that $p\in
(0,1)$ and $q>0.$  Then the following statements are equivalent:
\bnum \item[{\bf (\ref{p unit interval}a)}] $C_{\phi}$ is a compact
operator from ${\mathcal B}^p(U^n)$ to ${\mathcal B}^q(U^n)$.
\item[{\bf (\ref{p unit interval}b)}] $C_{\phi}$ is a compact operator from ${\mathcal
B}^p_0(U^n)$ to ${\mathcal B}^q(U^n)$. \item[{\bf (\ref{p unit
interval}c)}] $C_\phi$ is a bounded operator from ${\mathcal
B}^p(U^n)$ to ${\mathcal B}^q(U^n)$, and for all
$l\in\{1,2,\ldots,n\}$, we have that
\[
\lim_{|\phi_l(z)|\to 1}\sum\limits^n_{k=1}\left|\frac{\partial
\phi_{l}} {\partial
z_k}(z)\right|\frac{(1-|z_k|^2)^q}{(1-|\phi_l(z)|^2)^p}=0.
\]
\enum
\end{theorem}
Below are two consequences of Theorem \ref{p unit interval}:
\begin{corollary}\label{cor of t2and3}
Let $p,q>0$, and suppose that $\phi$ is a holomorphic self-map of
$U^n$, with $\|\phi_l\|_\infty<1$ and $\phi_l\in {\mathcal
B}^q(U^n)$ for each $l\in\{1,2,\ldots,n\}$.  Then $C_{\phi}$ is a
compact operator from ${\mathcal B}^p(U^n)$ and ${\mathcal
B}^p_0(U^n)$ to ${\mathcal B}^q(U^n)$.
\end{corollary}

\begin{corollary}\label{Stevo} Let $\phi$ be a holomorphic self-map
of $U^n$, with $p\in (0,1)$ and $q\geq 1.$  Then $C_{\phi}$ is a
compact operator from ${\mathcal B}^p(U^n)$ and ${\mathcal
B}^p_0(U^n)$ to ${\mathcal B}^q(U^n)$.
\end{corollary}

We postpone the statements and proofs of our corresponding
boundedness and compactness results for $C_{\phi}$ from ${\mathcal
B}^p_0(U^n)$ to ${\mathcal B}^q_0(U^n)$ until Section \ref{little
spaces}.  The following section is a gathering of the facts that
will be needed to prove the results that are stated in the present
section:

\section{Auxiliary Facts}\label{lemmas}
The following fact, which is used to prove Theorem
\ref{firsttheorem} and Lemmas \ref{hln1} and \ref{Macla}, does not
seem to appear in the literature; therefore, we give a proof for the
reader's convenience.

\begin{lemma}\label{pteval}
Let $z\in U^n$ and $p>0$.  Assume that $f\in{\mathcal B}^p(U^n).$
\bnum \item[{\bf (a)}] If $0<p<1,$ then
$$|f(z)|\leq \frac{n-p+1}{1-p}\|f\|_{\Bl}.$$

\item[{\bf (b)}] If $p=1$, then
$$|f(z)|\leq \frac{n\ln{2}+1}{n\ln{2}}\sum\limits^n_{k=1}\ln
\frac{2}{1-|z_k|^2}\|f\|_{\Bl}.$$

\item[{\bf (c)}] If $p>1,$ then
$$|f(z)|\leq \frac{2^{p-1}n+p-1}{n(p-1)}\sum\limits^n_{k=1}\frac{1}{(1-|z_k|^2)^{p-1}}\|f\|_{\Bl}.$$
\enum
\end{lemma}
\begin{proof}
Before verifying the individual parts of the lemma, we make a few
prerequisite observations.  Let $p>0$ and $z\in U^n$.  It follows
from the definition of $\|\cdot\|_{\Bl}$ that
\begin{equation}\label{conti}
|f(0)|\leq \|f\|_{\Bl},\qand\left|\frac{\partial f}{\partial
z_k}(z)\right|\leq \frac{\|f\|_{\Bl}}{(1-|z_k|^2)^p}
\end{equation}
for every $z\in U^n$ and $k\in\{1,2,\ldots,n\}$.  One also observes
that
\begin{small}
\begin{eqnarray*}f(z)-f(0)&=&\sum\limits^{n}_{k=1}f([0,z]_{n-k+1})-f([0,z]_{n-k})\\
&=&\sum\limits^{n}_{k=1}z_k\int^1_0\frac{\partial
f([0,(tz_k,z_k')]_{n-k+1})}{\partial z_k}dt.\end{eqnarray*}
\end{small}
It then follows from Inequalities (\ref{conti}) that
\begin{eqnarray}|f(z)|&\leq&
|f(0)|+\sum\limits^{n}_{k=1}|z_k|\int^1_0\frac{\|f\|_{\Bl}}{(1-|tz_k|^2)^p}dt\nonumber\\
&\leq&
\|f\|_{\Bl}+\|f\|_{\Bl}\sum\limits^{n}_{k=1}\int^{|z_k|}_0\frac{1}{(1-t^2)^p}dt.\label{start}
\end{eqnarray}
Let $k\in\{1,2,\ldots,n\}$.  Note that if $p=1$, then
\begin{equation}\label{p1}
\int^{|z_k|}_0\frac{1}{(1-t^2)^p}dt=\frac{1}{2}
\ln\frac{1+|z_k|}{1-|z_k|}\leq\frac{1}{2}\ln\frac{4}{1-|z_k|^2}.
\end{equation}
Also, if $p>0$ and $p\not=1$, then
\begin{eqnarray}&&\int^{|z_k|}_0\frac{1}{\left(1-t^2\right)^p}dt
\leq\int^{|z_k|}_0\frac{1}{(1-t)^p}dt=\frac{1-(1-|z_k|)^{-p+1}}{1-p}.\label{bunch}
\end{eqnarray}
We are now prepared to prove each part of the lemma:

{\bf (a):} By Inequality (\ref{bunch}), we have that \bege\label{ab}
\int^{|z_k|}_0\frac{1}{(1-t^2)^p}dt\leq\frac{1}{1-p}
\ende
for all $k\in\{1,2,\ldots,n\}$.  From Inequalities (\ref{start}) and
(\ref{ab}), it follows that
$$|f(z)|\leq\left(1+\frac{n}{1-p}\right)\|f\|_{\Bl}.$$  The statement in Part {\bf (a)}
is then obtained immediately.

{\bf (b):} Since
$$\ln\frac{4}{1-|z_k|^2}>\ln 4=2\ln 2$$
for each $k\in\{1,2,\ldots,n\},$ we have that \bege\label{inv}
1<\frac{1}{2n\ln 2}\sum\limits^n_{k=1}\ln\frac{4}{1-|z_k|^2}.
\ende
Combining Inequalities (\ref{start}), (\ref{p1}), and (\ref{inv}),
we then deduce that
$$|f(z)|\leq \left(\frac{1}{2}+\frac{1}{2n\ln
2}\right)\sum\limits^n_{k=1}\ln\frac{4}{1-|z_k|^2}\|f\|_{{\mathcal
B}^1},$$ from which the statement in Part {\bf (b)} of the lemma
follows.

{\bf (c):} By Inequality (\ref{bunch}), we have that for all
$k\in\{1,2,\ldots,n\}$,
\bege \int^{|z_k|}_0\frac{1}{(1-t^2)^p}dt\leq
\frac{1-(1-|z_k|)^{p-1}}{(p-1)(1-|z_k|)^{p-1}}
\leq\frac{2^{p-1}}{(p-1)(1-|z_k|^2)^{p-1}}.\label{en}
\ende
From Inequalities (\ref{start}) and (\ref{en}), we see that
\begin{eqnarray*}|f(z)|&\leq&\|f\|_{\Bl}+\frac{2^{p-1}}{p-1}\sum\limits^n_{k=1}
\frac{1}{(1-|z_k|^2)^{p-1}}\|f\|_{\Bl}\\
&\leq&
\left(\frac{1}{n}+\frac{2^{p-1}}{p-1}\right)\sum\limits^n_{k=1}
\frac{1}{(1-|z_k|^2)^{p-1}}\|f\|_{\Bl}.\end{eqnarray*} Therefore,
Part {\bf (c)} of the Lemma holds.
\end{proof}

The following lemma, which is based on a result due to Hardy and
Littlewood for $\alpha$-Lipschitz spaces defined using suprema on
$\partial U$ (rather than on $U$, as we need here), is essentially
Lemma 2 in \cite{mm1} and will be used several times in the proof
of Lemma \ref{hln1}.

\begin{lemma}\label{hln}
Let $0<p<1$.  Then ${\mathcal B}^{1-p}(U)={\Li}(U)$, and the norms
on these spaces are equivalent.
\end{lemma}
Lemma \ref{hln1} below and Theorem \ref{firsttheorem} together
imply Corollary \ref{lipresult}.  The lemma is also used to prove
Lemma \ref{usedonce?} and the necessity portion of Theorem \ref{p
unit interval}.
\begin{lemma}\label{hln1}
Lemma \ref{hln} holds with $U$ is replaced by $U^n$.
\end{lemma}

\begin{proof}
First, we show that ${\Li}(U^n)\subset{\mathcal B}^{1-p}(U^n)$ and
find a $C\geq 0$ such that for all $f\in {\Li}(U^n)$,
$||f||_{{\mathcal B}^{1-p}}\leq C ||f||_{\Li}$.  Let
$f\in{\Li}(U^n)$, $j\in\{1,2,\ldots,n\}$, and $z,w\in U^n$, so that
$|f(z)-f(w)|\leq ||f||_{\Li}|z-w|^p$. Since $z$ and $w$ were chosen
arbitrarily, we can replace $w$ by $(w_j,z_j')$ in this inequality.
We then obtain that
\[|f(z)-f(w_j,z_j')|\leq ||f||_{\Li}|z_j-w_j|^p.
\]
\noindent Define $f_{z'_j}:U\rightarrow\CC$ by
\[
f_{z'_j}(\zeta)=f(\zeta,z_j'),
\]
and note that $f_{z'_j}\in{\Li}(U)$.  This fact and the estimates
\begin{eqnarray*}
|f_{z_j'}(0)|&\leq &|f_{z_j'}(0)-f(0)|+|f(0)|\\
&\leq &||f||_{\Li}|z'_j|^\alpha+||f||_{\Li}\\
&\leq &(1+n^{\frac{\alpha}{2}})||f||_{\Li}
\end{eqnarray*}
together imply that $||f_{z_j'}||_{{\Li}}\leq
(2+n^{\frac{\alpha}{2}})||f||_{\Li}$. By Lemma \ref{hln}, we then
have that $f_{z'_j}\in{\mathcal B}^{1-p}(U)$ and the stronger fact
that there exists a $C'\geq 0$ such that for all
$j\in\{1,2,\ldots,n\}$, $z\in U^n$, and $f\in\Li(U^n)$,
\[
||f_{z'_j}||_{{\mathcal B}^{1-p}}\leq C'||f||_{\Li}.
\]
We rewrite the inequality above in the form
\[
\sup_{w\in U }\left|\frac{df_{z'_j}}{d\zeta}
(\zeta)\right|(1-|\zeta|^2)^{1-p}\leq C'||f||_{\Li}.
\]
Since $z\in U^n$ was chosen arbitrarily, we have that for all
$f\in{\Li}(U^n)$ and $j\in\{1,2,\ldots,n\}$,
\[
\sup_{z\in U^n}\left|\frac{\partial f_{z'_j}}{\partial
z_j}(z_j)\right|(1-|z_j|^2)^{1-p}\leq C'||f||_{\Li}.
\]
Hence, we have that
\begin{equation}\label{lipo}
{\Li}(U^n)\subset{\mathcal B}^{1-p}(U^n),
\end{equation}
and $||f||_{{\mathcal B}^{1-p}}\leq C||f||_{\Li}$ for all $f\in
{\Li}(U^n)$, where $C=1+C'n$.

Next, we show that ${\mathcal B}^{1-p}(U^n)\subset{\Li}(U^n)$ and
that there is a $C''\geq 0$ such that for all $f\in {\mathcal
B}^{1-p}(U^n)$,
\begin{equation}\label{lipcon}
||f||_{\Li}\leq C''||f||_{{\mathcal B}_{1-p}}.
\end{equation}
\noindent Of course, by the set inclusion (\ref{lipo}), it will
follow that the above inequality holds for all $f\in{\Li}(U^n)$ as
well, and the proof of the theorem will be complete. If
$j\in\{1,2,\ldots,n\}$, $f\in{\mathcal B}^{1-p}(U^n)$, and $z\in
U^n$, then
\[
\left|\frac{\partial f_{z'_j}}{\partial
z_j}(z_j)\right|(1-|z_j|^2)^{1-p}\leq ||f||_{{\mathcal B}^{1-p}}.
\]
\noindent It follows that
\[
\sup_{z,w\in U^n}\left|\frac{\partial
f_{([z,w]_{n-j+1})'_j}}{\partial z_j}(z_j)\right|(1-|z_j|^2)^{1-p}
\leq ||f||_{{\mathcal B}^{1-p}}.
\]
Hence, for any choice of $z,w\in U^n$, we have that
\bege\label{529}
\sup_{\lambda\in U}\left|\frac{d
f_{([z,w]_{n-j+1})'_j}}{d\lambda}(\lambda)\right|(1-|\lambda|^2)^{1-p}\leq
||f||_{{\mathcal B}^{1-p}},
\ende
\noindent which implies that $ f_{([z,w]_{n-j+1})'_j}\in{\mathcal
B}^{1-p}(U)$ for all $z,w\in U^n$.  Now Lemma \ref{pteval} implies
that
\[
|f_{([z,w]_{n-j+1})'_j}(0)|\leq \frac{n+p}{p}||f||_{{\mathcal
B}^{1-p}}.
\]
Combining this fact with Inequality (\ref{529}), we obtain that
\[
|| f_{([z,w]_{n-j+1})'_j}||_{{\mathcal B}^{1-p}}\leq
\frac{n+2p}{p}||f||_{{\mathcal B}^{1-p}}.
\]
\noindent It follows by a further application of Lemma \ref{hln}
that there is a $C'''\geq 0$ such that for all $z,w\in U^n$ and
$f\in {\mathcal B}^{1-p}(U^n)$,
\[
|| f_{([z,w]_{n-j+1})'_j}||_{{\mathcal L}_p}\leq
C'''||f||_{{\mathcal B}^{1-p}}.
\]
\noindent  We now have that for all $\lambda,\nu\in U$, $z,w\in
U^n$, and $f\in {\mathcal B}^{1-p}(U^n)$,
\[
|f_{([z,w]_{n-j+1})_j'}(\lambda)-f_{([z,w]_{n-j+1})_j'}(\nu)|\leq
C'''||f||_{{\mathcal B}^{1-p}}|\lambda-\nu|^p.
\]
\noindent Putting $\lambda=z_j$ and $\nu=w_j$, we obtain that
\[
|f([z,w]_{n-j})-f([z,w]_{n-j+1})|\leq C'''||f||_{{\mathcal
B}^{1-p}}|z_j-w_j|^p.
\]
\noindent By the telescoping inequality
\[
|f(z)-f(w)|\leq \sum_{j=1}^{n}|f([z,w]_{n-j})-f([z,w]_{n-j+1})|,
\]
\noindent we have that for all $z,w\in U^n$ and $f\in{\mathcal
B}^{1-p}(U^n)$,
\begin{eqnarray*}
|f(z)-f(w)|&\leq &C'''||f||_{{\mathcal B}^{1-p}}\sum_{j=1}^n|z_j-w_j|^p\\
&\leq & nC''||f||_{{\mathcal B}^{1-p}}|z-w|^p.
\end{eqnarray*}
\noindent Hence, ${\mathcal B}^{1-p}(U^n)\subset {\Li}(U^n)$ and
therefore equals ${\Li}(U^n)$ by the set inclusion (\ref{lipo});
furthermore, the above inequality implies that there exists a
$C''=nC''\geq 0$ such that for all $f\in {\mathcal
B}^{1-p}(U^n)={\Li}(U^n)$,
\[||f||_{\Li}\leq C''||f||_{{\mathcal B}^{1-p}}.
\]
Combining these facts with Inequality (\ref{lipcon}), we obtain the
statement of the lemma.
\end{proof}
The following elementary consequence of
\cite[~Thm.~5.2,~p.~302]{Dug} will be used in the proofs of Lemma
\ref{usedonce?} and necessity in Theorem \ref{p unit interval}:
\begin{lemma}\label{Cauchy}
Suppose that $D\subset\CC^n$ is a bounded domain, and let $X$ be a
Banach space with norm $||\cdot||_X$.  Let $\alpha\in(0,1]$, and
assume that $f:D\rightarrow X$ has the property that there is a
$C\geq 0$ such that for all $z,w\in D$,
\[
||f(z)-f(w)||_X\leq C|z-w|^\alpha.
\]
Then $f$ has a unique, uniformly continuous extension $g:{\overline
D}\rightarrow X$.
\end{lemma}
The following fact is used in the proofs of Lemmas \ref{precurs506}
and \ref{littlelemma}, Theorem \ref{p unit interval} and Theorem
\ref{secondtheorem}, Part {\bf (b)}.

\begin{lemma}\label{Macla} Let $\phi:U^n\rightarrow U^n$ be holomorphic, and suppose that
$p,q>0$.  Then $C_{\phi}$ is a compact operator from ${\mathcal
B}^p(U^n)$ (respectively, $ {\mathcal B}^p_0(U^n)$) to ${\mathcal
B}^q(U^n)$ if and only if for any bounded sequence
$(f_j)_{j\in\NN}$ in ${\mathcal B}^p(U^n)$ (respectively,
${\mathcal B}^p_0(U^n)$) such that $f_j\to0$ uniformly on compacta
in $U^n$ as $j\rightarrow\infty$, we have that
$\|C_{\phi}f_j\|_{{\mathcal B}^q}\to 0$ as $j\to\infty.$
\end{lemma}

\begin{proof} $\Rightarrow$: Assume that $C_\phi$ is compact, and let $(f_j)_{j\in\NN}$ be a
bounded sequence in ${\mathcal B}^p(U^n)$, (respectively,
${\mathcal B}^p_0(U^n)$) with $f_j\to 0$ uniformly on compacta in
$U^n.$ Suppose to the contrary that there is a subsequence
$(f_{j_m})_{m\in\NN}$ and a $\d>0$ such that
$\|C_\phi(f_{j_m})\|_{{\mathcal B}^q}\geq \d$ for all $m\in
{\NN}.$ Since $C_\phi$ is compact, we have that
$(C_\phi(f_{j_m}]))_{m\in\NN}=(f_{j_m}\circ \phi)_{m\in\NN}$ has a
further subsequence $(f_{j_{m_l}}\circ\phi)_{l\in\NN}$ such that
for some $g\in {\mathcal B}^q(U^n)$,
\begin{equation}\label{conver}
\lim_{l\rightarrow\infty}||f_{j_{m_l}}\circ\phi-g||_{{\mathcal
B}^q}=0.
\end{equation}

Lemma \ref{pteval} implies that for any compact $K\subset U^n$ there
is a $C_K\geq 0$ such that for all $l\in\NN$ and $z\in K$, the
following inequality holds:
\begin{equation}\label{domin}
|f_{j_{m_l}}(\phi(z))-g(z)|\leq C_K\|f_{j_{m_l}}\circ
\phi-g\|_{{\mathcal B}^q}.
\end{equation}
Therefore, by Equation (\ref{conver}), $f_{j_{m_l}}\circ\phi-g\to
0$ uniformly on compacta in $U^n$.  Now $f_{j_{m_l}}(\phi(z))\to
0$ as $l\to\infty$ for each $z\in U^n$, by hypothesis and the fact
that $\{\phi(z)\}$ is compact for each $z\in U^n$.  By Inequality
(\ref{domin}), we must then have that $g=0$; therefore, by
Equation (\ref{conver}), it follows that
$$\lim_{l\to\infty}\|C_\phi(f_{j_{m_l}})\|_{{\mathcal B}^q}=0.$$
We have now obtained a contradiction.

$\Leftarrow$:  Assume that $(g_j)_{j\in\NN}$ is a sequence in
${\mathcal B}^p$ such that $||g_j||_{\Bl}\leq M$ for all
$j\in\NN$.
 Lemma \ref{pteval} implies that $(g_j)_{j\in\NN}$ is uniformly bounded on
compacta in $U^n$ and consequently normal by Montel's theorem. We
then have that there is a subsequence $(g_{j_m})_{m\in\NN}$ of
$(g_j)_{j\in\NN}$ that converges uniformly on compacta in $U^n$ to
some $g\in {\mathcal H}(U^n).$ It follows that ${\pt g}_{j_m}/{\pt
z_l}\to {\pt g}/{\pt z_l}$ uniformly on compacta in $U^n$ for each
$l\in\{1,2,\ldots,n\}.$ Thus $g\in {\mathcal B}^p(U^n)$, and
$\|g_{j_m}-g\|_{\Bl}\leq M+||g||_{\Bl}<\infty$ and converges to $0$
on compacta in $U^n.$  By hypothesis, we then have that
$g_{j_m}\circ\phi\to g\circ\phi$ in ${\mathcal B}^q(U^n).$
Therefore, $C_\phi$ is a compact operator from ${\mathcal
B}^q(U^n)$.
\end{proof}
The following observation is used in the proof of Theorem \ref{p
unit interval}:
\begin{lemma}\label{usedonce?}
Let $p\in (0,1)$, and assume that $(f_j)_{j\in\NN}$ is a bounded sequence in ${\mathcal B}^p(U^n)$ such that
$f_j\rightarrow 0$ uniformly on compacta in $U^n.$ Then the
continuous extensions $f_j\rightarrow 0$ uniformly on $\ol{U^n}.$
\end{lemma}
\begin{proof}
By Lemma \ref{hln1}, $(f_j)_{j\in\NN}\subset{\mathcal
L}_{1-p}(U^n)$, and there is a $C\geq 0$ such that for all $j\in
\NN$, $||f_j||_{{\mathcal L}_{1-p}}\leq C$. The fact that each $f_j$
extends continuously to $\ol{U^n}$ then follows from Lemma
\ref{Cauchy}. If $C=0$, then each $f_j$ is zero and trivially,
$f_j\rightarrow 0$ uniformly on $\ol{U^n}$. Therefore, we can now
assume that $C>0$. Let $\ve\in(0,1)$, and put
\[
\delta=\min\left\{\frac{\ve}{2},\,\,\,\frac{1}{n}\left(\frac{\ve}{2C}\right)^{\frac{1}{1-p}}\right\}.
\]
Since $(1-\delta)\ol{U^n}$ is compact, then by hypothesis there is
an $N\in\NN$ such that for all $z\in (1-\delta)\ol{U^n}$ and
$j\in\NN$ with $j\geq N$, $|f_j(z)|<\ve/2$.  Hence, if we can show
that $|f_j(z)|<\ve$ for all $z\in \ol{U^n}-(1-\delta)\ol{U^n}$ and
$j\geq N$, then the proof will be complete.  Letting
$w_z=[(1-\delta)/|z|]z$ for each $z\in \ol{U^n}-(1-\delta)\ol{U^n}$,
we see that $|z-w_z|<n\delta$ for all $z\in\ol{U^n}$.  Since
$||f_j||_{\Li}\leq C$ for all $j\in\NN$, we must have that for all
such $j\geq N$ and $z\in \ol{U^n}-(1-\delta)\ol{U^n}$,
\[
|f_j(z)|\leq |f_j(w_z)|+C|z-w_z|^{1-p}\leq
\frac{\ve}{2}+C|n\delta|^{1-p} < \ve.
\]
\end{proof}
We now obtain norm estimates for a certain class of testing
functions, which will be used in the proofs of Theorems
\ref{secondtheorem}, \ref{p unit interval}, and Lemma
\ref{precurs3}.
\begin{lemma}\label{test} Let $p>0$.  Then ${\mathcal
T}^p\subset{\mathcal B}_0^p(U^n)$.  Furthermore, there is a $Q\geq
0$ such that $||\nu||_{{\mathcal B}^p}\leq Q$ for all
$\nu\in{\mathcal T}^p$.
\end{lemma}

{\it Proof.} Let $k,l\in\{1,2,\ldots,n\}$, and pick $w\in U$. Then
we have that \bege\label{par} \frac{\pt f^{(l)}_w}{\pt
z_k}=0{\text{ on }}U^n\,\,\,{\text{  for }}k\not=l,
\ende
and \bege\label{par1} \frac{\pt f^{(l)}_w}{\pt z_l}(z)=\frac{
1}{(1-\bar wz_l)^p}\,\,\,{\text{  for all }}z\in U^n.
\ende
Therefore, for all $z\in U^n$, we have that
$$|f^{(l)}_w(0)|+\sum\limits^n_{k=1}
\left|\frac{\partial f^{(l)}_w} {\partial
z_k}(z)\right|(1-|z_k|^2)^p=\frac{(1-|z_l|^2)^p}{|1-\bar
wz_l|^p}\leq (1+|z_l|)^p\leq 2^p.$$ It follows that ${\mathcal
T}^p_1\subset {\mathcal B}^p(U^n)$, and
\begin{equation}\label{nu1}
||\nu||_{{\mathcal B}^p}\leq 2^p{\text{ for all }}\nu\in {\mathcal
T}^p_1.
\end{equation}

We now prove that ${\mathcal T}^p_1\subset{\mathcal B}_0^p(U^n)$.
Let $l\in\{1,2,\ldots,n\}$.  Since we have that
$$(1-\bar w t)^{-p}=\sum\limits^{+\infty}_{j=0}\frac{p(p+1)\cdots
(p+j-1)}{j!}(\bar w)^jt^j,$$
Equation (\ref{par1}) implies that
$$f_w^{(l)}(z)=\sum\limits^{+\infty}_{j=0}\frac{p(p+1)\cdots
(p+j-1)}{j!}(\bar w)^j\int^{z_l}_0 t^jdt.$$ For $m\in\WW$, let
$P^{(l,1)}_{m,w}:\CC^n\rightarrow\CC$ be the polynomial function
defined by
$$P_{m,w}^{(l,1)}(z)=\sum\limits^{m}_{j=0}\frac{p(p+1)\cdots (p+j-1)}{j!}(\bar w)^j\int^{z_l}_0t^jdt.$$
We then have by Equations (\ref{par}) and (\ref{par1}) that for
$m\in\WW$ and $z\in U^n$,
\begin{eqnarray*}&&|f_w^{(l)}(0)-P_{m,w}^{(l,1)}(0)|+\sum_{k=1}^n\left|\frac{\partial [f_w^{(l)}-P_{m,w}^{(l,1)}]}{\partial
z_k}(z)\right|(1-|z_k|^2)^p \\
&=&|f_w^{(l)}(0)-P_{m,w}^{(l,1)}(0)|+\left|\frac{\partial
[f_w^{(l)}-P_{m,w}^{(l,1)}]}{\partial
z_l}(z)\right|(1-|z_l|^2)^p \\
&\leq& \sum\limits^{+\infty}_{j=m+1}\frac{p(p+1)\cdots
(p+j-1)}{j!}|w|^j\to 0 \;\mbox{as}\;\; m\to \infty.
\end{eqnarray*}
It follows that $f_w^{(l)}\in {\mathcal B}^p_0(U^n)$, thus showing
that
\begin{equation}\label{t1}
{\mathcal T}_1^p\subset{\mathcal B}^p_0(U^n).
\end{equation}

Next, consider ${\mathcal T}^p_2$.  Let $w\in U$ and
$l\in\{1,2,\ldots,n\}$.  Then both of the following equations hold
for all $z\in U^n$:
\begin{equation}\label{gzero}
\frac{\partial g^{(l)}_w}{\partial z_{k}}(z)=0\quad ({\text{for
}}k\neq l);\quad \frac{\partial g^{(l)}_w}{\partial z_l}(z)=
p\bar{w}\frac{1-|w|^2}{(1-z_l\bar{w})^{p+1}}.
\end{equation}
By an elementary estimate, it follows that for all $z\in U^n$,
$$|g^{(l)}_w(0)|+\sup\limits_{z\in U^n}\sum\limits^n_{k=1}
(1-|z_k|^2)^p\left|\frac{\partial g^{(l)}_w}{\partial
z_k}(z)\right|\leq 1+p2^{p+1}.$$ Therefore, ${\mathcal
T}^p_2\subset{\mathcal B}^p(U^n)$, and we have that
\begin{equation}\label{nu2}
||\nu||_{{\mathcal B}^p}\leq 1+p2^{p+1}{\text{ for all
}}\nu\in{\mathcal T}^p_2.
\end{equation}
We now show that \bege\label{t2} {\mathcal T}^p_2\subset{\mathcal
B}^p_0(U^n).
\ende
We can write
$$\frac{1}{(1-z_l\bar{w})^p}=\sum_{j=0}^\infty\frac{p(p+1)\cdots(p+j-1)}{j!}
\left(\bar{w}z_l\right)^j.$$ For each $m\in\WW$, define the
polynomial function $P_{m}^{(l,2)}:\CC^n\rightarrow\CC$ by
$$P_{m}^{(l,2)}(z)=(1-|w|^2)\sum_{j=0}^m\frac{p(p+1)\cdots(p+j-1)}{j!}\left(\bar{w}z_l\right)^j.$$
Arguing as in the case of $P_{m,w}^{(l,1)}$ above, we obtain that
$||g^{(l)}_w-P_{m,w}^{(l,2)}||_{\Bl}\rightarrow 0$ as
$m\rightarrow\infty$.  The desired set inclusion (\ref{t2}) follows.

Next, we consider ${\mathcal T}^p_3$.  Suppose that
$l\in\{2,3,\ldots,n\}$, and let $w\in U$.
 Since
 \bege\label{h0}\frac{\partial h_w^{(l)}}{\partial
z_{k}}=0\,\,\,{\text{ for }}k\neq1,l,
\ende
\bege\label{hone} \frac{\partial h_w^{(l)}}{\partial
z_1}=\left(\frac{1-|w|^2}{1-z_l\bar{w}}\right)^p,{\text{ and }}
\ende
\bege\label{helle} \frac{\partial h_w^{(l)}}{\partial
z_l}=p(z_1+2)\bar{w}\frac{(1-|w|^2)^p}{(1-z_l\bar{w})^{p+1}},
\ende
we have that
$$|h_w^{(l)}(0)|+\sup\limits_{z\in U^n}\sum\limits^n_{k=1}(1-|z_k|^2)^p
\left|\frac{\partial h_w^{(l)}}{\partial z_k}\right|\leq
2+2^p+3p2^{p+1}<\infty.$$ The above estimate implies that
${\mathcal T}^p_3\subset{\mathcal B}^p(U^{n})$ and that
\bege\label{nu3} ||\nu||_{{\mathcal B}^p}\leq
2+2^p+3p2^{p+1}\,\,\,{\text{ for all }}\nu\in{\mathcal T}^p_3.
\ende
By a similar argument to the one that we used for $g_w^{(l)}$,
$h_w^{(l)}$ is the ${\mathcal B}^p$-norm limit (as
$m\rightarrow\infty$) of the functions
$P_{m,w}^{(l,3)}:\CC^n\rightarrow\CC$ given by
$$P_{m,w}^{(l,3)}(z)=(z_1+2)(1-|w|^2)^p\sum_{j=0}^m
\frac{p(p+1)\cdots(p+j-1)}{j!}\bar{w}^jz_l^j,$$ and we thus obtain
that \bege\label{t3} {\mathcal T}^p_3\subset{\mathcal B}^p_0(U^n).
\ende
The first statement in the lemma now follows from the set
inclusions (\ref{t1}), (\ref{t2}), and (\ref{t3}).  The second
statement in the lemma is obtained from Inequalities (\ref{nu1}),
(\ref{nu2}), and (\ref{nu3}).

\medskip

The following observation is used in the proofs of Theorems
\ref{secondtheorem}, \ref{p unit interval}, \ref{almostmissed}, and
Corollary \ref{firstiff}.
\begin{lemma}\label{phicomp} If $C_{\phi}$ maps
${\mathcal B}^p(U^n)$ (respectively, ${\mathcal B}^p_0(U^n)$) to
${\mathcal B}^q(U^n)$ (respectively, ${\mathcal B}^q_0(U^n)$), then
$\phi^\gamma\in{\mathcal B}^q(U^n)$ (respectively, ${\mathcal
B}^p_0(U^n)$) for each $n$-multi-index $\gamma$.
\end{lemma}
\begin{proof}
The polynomial function $p_\gamma:\CC^n\rightarrow\CC$ defined by
$p_\gamma(z)=z^\gamma$ is of course in ${\mathcal B}^p_0(U^n)$, so
by assumption, $C_\phi(p_\gamma)=\phi^\gamma\in {\mathcal B}^q(U^n)$
(respectively, ${\mathcal B}^q_0(U^n)$).
\end{proof}
Proposition \ref{precurs506} below is used in the proof of Theorem
\ref{secondtheorem}, Part {\bf (a)}:
\begin{proposition}\label{precurs506}
Let $p,q>0$, and assume that $\phi:U^n\rightarrow U^n$ is a
holomorphic map such that $\phi_k\in {\mathcal B}^q(U^n)$ for each
$k\in\{1,2,\ldots,n\}$ and Condition $(\ref{firstcom})$ holds.
Then $C_\phi$ is a compact operator from ${\mathcal B}^p(U^n)$ and
${\mathcal B}^p_0(U^n)$ to ${\mathcal B}^q(U^n)$.
\end{proposition}
\begin{proof}
Let $(f_j)_{j\in\NN}$ be a sequence in ${\mathcal B}^p(U^n)$
(respectively, ${\mathcal B}^p_0(U^n)$) such that $f_j\rightarrow 0$
uniformly on compacta in $U^n$ and
\begin{equation}\label{bound5}
\left\|f_{j}\right\|_{\Bl}\leq C\,\,\,{\text{ for all }}j\in\NN.
\end{equation}
By Lemma \ref{Macla}, it suffices to show that
\begin{equation}\label{finish2}
\lim_{j\rightarrow\infty}\left\|C_{\phi}f_{j}\right\|_{{\mathcal
B}^q}=0.
\end{equation}
If $C=0$, then $(f_j)_{j\in\NN}$ and $(C_\phi f_j)_{j\in\NN}$ are
both zero sequences, and we trivially obtain Equation
(\ref{finish2}). Therefore, we can assume that $C>0$.  If
$||\phi_m||_{{\mathcal B}^q}=0$ for all $m\in\{1,2,\ldots,n\}$, then
$\phi=0$ and $C_\phi$ has finite rank. Therefore, we can also assume
that $||\phi_m||_{{\mathcal B}^q}>0$ for some
$m\in\{1,2,\ldots,n\}$.

 Let $\varepsilon>0$.  By the assumed Condition (\ref{firstcom}), there is an $r\in(0,1)$
such that for all $z\in U^n$ satisfying $d(\phi(z),\partial
U^n)<r$,
\begin{equation}\label{peper}
\sum\limits^n_{k,l=1}\left|\frac{\partial \phi_l}{\partial
z_k}(z)\right| \frac{(1-|z_k|^2)^q}{(1-|\phi_l(z)|^2)^p}
<\frac{\varepsilon}{2C}.\end{equation} Subsequent applications of
the chain rule, and Inequalities (\ref{bound5}) and (\ref{peper})
yield that for all such $z$ and all $j\in\NN$,
\begin{eqnarray}
&&\sum\limits^n_{k=1}\left|\frac{\partial\left(C_{\phi}f_j\right)}
{\partial z_k}(z)\right| (1-|z_k|^2)^q\nonumber\\
&\leq& \sum\limits^n_{l=1}\left|\frac{\partial f}{\partial
\zeta_l}[\phi(z)]\right|\left(1-|\phi_l(z)|^2\right)^p\sum\limits^n_{k=1}\left|\frac{\partial
\phi_{l}} {\partial z_k}(z)\right|\frac{\left(1-|z_k|^2\right)^q}
{\left(1-|\phi_l(z)|^2\right)^p}\nonumber\\
&\leq&C\sum\limits^n_{k,l=1}\left|\frac{\partial \phi_l}{\partial
z_k}(z)\right|\frac{(1-|z_k|^2)^q}{(1-|\phi_l(z)|^2)^p}\nonumber\\
&<&\frac{\varepsilon}{2}\label{bound87}.
\end{eqnarray}
We now obtain the same estimate in the case that $d(\phi(z),\partial
U^n)\geq r.$ Let $E_r=\{w\,| \, d(w,\partial U^n)\geq r\}$.  $E_r$
is compact; therefore, by hypothesis, $(f_j)_{j\in\NN}$ and in turn
the sequences of partial derivatives $\left(\partial f_j/\partial
z_l\right)_{j\in\NN}$ for $l\in\{1,2,\ldots,n\}$ converge to $0$
uniformly on $E_r$. Consequently, there is an $N_1\in\NN$ such that
for all $w\in E_r$ and $j\in\NN$ with $j\geq N_1$,
\[
\sum_{l=1}^n\left|\frac{\partial f_j}{\partial
\zeta_l}(w)\right|<\frac{\varepsilon}{2\sum_{m=1}^n||\phi_m||_{{\mathcal
B}^q}}.
\]
It follows that,
\begin{small}
\begin{eqnarray}
\sum\limits^n_{k=1}\left|\frac{\partial\left(C_{\phi}f_j\right)}
{\partial z_k}(z)\right| (1-|z_k|^2)^q&\leq&
\sum\limits^n_{l=1}\left|\frac{\partial f_j}{\partial
\zeta_l}(\phi(z))\right|\sum\limits^n_{k=1}\left|\frac{\partial
\phi_{l}}
{\partial z_k}(z)\right|{\left(1-|z_k|^2\right)^q}\nonumber\\
&\leq &\sum\limits^n_{l=1}\left|\frac{\partial f_j}{\partial
\zeta_l}(\phi(z))\right|||\phi_l||_{{\mathcal B}^q(U^n)}\nonumber\\
&\leq &\sum\limits^n_{l=1}\sup_{w\in E_r}\left|\frac{\partial
f_j}{\partial
\zeta_l}(w)\right|||\phi_l||_{{\mathcal B}^q(U^n)}\nonumber\\
&<&\frac{\varepsilon}{2}.\label{bound88}
\end{eqnarray}
\end{small}
Since $\{\phi(0)\}$ is compact, it must be the case that
$f_j(\phi(0))\rightarrow 0$ as $j\rightarrow\infty$, so there is
an $N_2\in\NN$ such that for all $j\in\NN$ with $j\geq N_2$,
\begin{equation}\label{pointwi}
|f_j(\phi(0))|<\frac{\varepsilon}{2}.
\end{equation}
The proposition then follows if we let $N=\max(N_1,N_2)$ and apply
Inequalities (\ref{bound87}), (\ref{bound88}), and
(\ref{pointwi}).\msk

\end{proof}

The following lemma is used twice in the proof of Theorem
\ref{secondtheorem}:
\begin{lemma}\label{precurs3}
Suppose that $\phi:U^n\rightarrow U^n$ is holomorphic, and assume
that $p,q>0$.  Let $C_\phi$ be bounded from ${\mathcal B}^p(U^n)$
(respectively, ${\mathcal B}^p_0(U^n)$) to ${\mathcal B}^q(U^n)$.
Suppose that there is a sequence $(z^{[j]})_{j\in\NN}$ in $U^n$ such
that $\om^{[j]}:=\phi(z^{[j]})\to \pt U^n$ as $j\to\infty$, and an
$\ve>0$ such that
\begin{equation}\label{attaboy}
\sum\limits^n_{k,l=1}\left|\frac{\partial \phi_{l}}{\partial
z_k}(z^{[j]})
\right|\frac{(1-|z^{[j]}_k|^2)^q}{(1-|\phi_l(z^{[j]})|^2)^p} \geq
\ve,\end{equation} for all $j\in\NN.$  If $l\in\{1,2,\ldots,n\}$
is such that $(\om^{[j]}_l)_{j\in\NN}$ has a subsequence
$(\om^{[j']}_l)_{j\in\NN}$ that converges to some $\rho_l$ with
$|\rho_l|=1$, then there is a sequence $(f_j)_{j\in\NN}$ in
${\mathcal B}^p(U^n)$ with the following properties:

\bnum \item[{\bf (i)}] $(||f_j||_{\mathcal B^p})_{j\in\NN}$ is
bounded, and $f_j\in{\mathcal B}^p_0(U^n)$ for all $j\in\NN$;

\item[{\bf (ii)}] $f_j\rightarrow 0$ uniformly on compacta in
$U^n;$

\item[{\bf (iii)}] $\|C_\phi f_j\|_{{\mathcal B}^q}\not\to 0,$ as
$j\to\infty.$ \enum
\end{lemma}
\begin{proof}
Since $C_\phi$ is bounded from ${\mathcal B}^p(U^n)$
(respectively, ${\mathcal B}^p_0(U^n)$) to ${\mathcal B}^q(U^n)$,
Condition {\bf (\ref{firsttheorem}c)} holds. It follows from
Inequality (\ref{attaboy}) that there is a subsequence
$(z^{[j']})_{j\in\NN}$ of $(z^{[j]})_{j\in\NN}$ and an
$\varepsilon'>0$ such that for some fixed
$k_0,l\in\{1,2,\ldots,n\}$,
$$\lim_{j\rightarrow\infty}
\left|\frac{\partial \phi_{l}}{\partial z_{k_0}}(z^{[j']})
\right|\frac{(1-|z^{[j']}_{k_0}|^2)^q}{(1-|\phi_l(z^{[j']})|^2)^p}=\varepsilon'.$$
Define $f_j:\ol{U^n}\rightarrow\CC$ for each $j\in\NN$ by
\[f_{j}(z)=g_{\om_l^{[j']}}^{(l)}(z).\]
Noticing that $f_j\in{\mathcal T}^p_2$ for all $j\in\NN$, we can
apply Lemma \ref{test} to obtain that $(f_j)_{j\in\NN}$ is a bounded
sequence in ${\mathcal B}^p(U^n)$ and also lies in ${\mathcal
B}^p_0(U^n)$. Therefore, $(f_j)_{j\in\NN}$ has the property ({\bf
i}). An elementary estimate shows that this sequence also satisfies
property ({\bf ii}).

To prove that property ({\bf iii}) holds, observe that
\begin{small}
\[
\sum\limits^n_{k=1}\frac{\partial [C_\phi(f_j)]}{\partial
z_k}(z^{[j']})(1-|z^{[j']}_k|^2)^q=\sum\limits^n_{k=1}(1-|z^{[j']}_k|^2)^q
\left|\frac{\partial f_j}{\partial
\zeta_l}(\phi(z^{[j']}))\frac{\partial \phi_l}{\partial
z_k}(z^{[j']})\right|
\]
\end{small}
\begin{eqnarray*}
=&p|\om^{[j']}_l|\sum\limits^n_{k=1}\frac{(1-|z^{[j']}_k|^2)^q}{(1-|\om^{[j']}_l|^2)^p}
\left|\frac{\partial \phi_l}{\partial z_k}(z^{[j']})\right|\\
\geq&p|\om^{[j']}_l|\left|\frac{\partial \phi_l}{\partial
z_{k_0}}(z^{[j']})
\right|\frac{(1-|z^{[j']}_{k_0}|^2)^q}{(1-|\phi_l(z^{[j']})|^2)^p}\\
\to& p\ve'>0\;\; \mbox{as}\;\; j\to\infty.
\end{eqnarray*}
The proof of the Lemma \ref{precurs3} is now complete.
\end{proof}
The following result, which we will use to prove Corollary
\ref{Stevo}, appears in \cite{zs3} and in more general form in
\cite{t}. Recall that a domain in $\CC^n$ is called {\em
homogeneous} \cite{t} if and only if it has a transitive
automorphism group.
\begin{theorem}[Generalized Schwarz Lemma]\label{Schwarzlemmapolydisk} Let $D_1$, $D_2$ be bounded,
homogeneous domains in $\CC^n$, and let $H^{D_1}$ and $H^{D_2}$
denote their respective Bergman metrics. Let $\phi:D_1\rightarrow
D_2$ be holomorphic.  Then there exists a $C>0$ such that for all
$z\in D_1$ and $u\in\CC^n$,
\[
H_{\phi(z)}^{D_2}\big(J_\phi(z)u,{\overline{ J_\phi(z)}}\big)\leq
CH_z^{D_1}(u,{\bar u}).
\]
\end{theorem}
Theorem \ref{Schwarzlemmapolydisk} above concludes the list of
prerequisite facts that will be used to prove our main results.
\section{Proofs of Main Results}\label{proofs}

In this section, we prove Theorems \ref{firsttheorem},
\ref{secondtheorem}, \ref{p unit interval} and their
corollaries.\msk

{\it Proof of Theorem \ref{firsttheorem}.}  First assume that
Condition {\bf (\ref{firsttheorem}c)} holds.  By Lemma \ref{pteval},
there is a $C\geq 0$ such that for all $f\in {\mathcal B}^p(U^n)$,
\bege\label{ptphi} |f(\phi(0))|\leq C\|f\|_{\Bl}.
\ende
The following chain of inequalities also holds for each $z\in U^n$:
\begin{eqnarray*}
&&\sum\limits^n_{k=1}\left|\frac{\partial \left[C_{\phi}f(z)\right]}
{\partial z_k}\right| (1-|z_k|^2)^q\\
&\leq&\sum\limits^n_{l=1}\left|\frac{\partial f}{\partial
\zeta_l}[\phi(z)]\right|\left(1-|\phi_l(z)|^2\right)^p\sum\limits^n_{k=1}\left|\frac{\partial
\phi_{l}} {\partial z_k}(z)\right|\frac{\left(1-|z_k|^2\right)^q}
{\left(1-|\phi_l(z)|^2\right)^p}\\
&\leq &||f||_{\Bl}\sum\limits^n_{k,l=1}\left|\frac{\partial
\phi_{l}} {\partial z_k}(z)\right|\frac{\left(1-|z_k|^2\right)^q}
{\left(1-|\phi_l(z)|^2\right)^p}\\
&\leq& M||f||_{\Bl}.
\end{eqnarray*}
Conditions {\bf (\ref{firsttheorem}a)} and {\bf
(\ref{firsttheorem}b)} then follow from simultaneous applications of
Inequality (\ref{ptphi}) and the inequality above.

Next, let $l\in\{1,2,\ldots,n\}$, and assume that Condition {\bf
(\ref{firsttheorem}a)} (respectively, {\bf (\ref{firsttheorem}b)})
holds; i.e., suppose that there is a $C\geq 0$ such that for all
$f\in {\mathcal B}^p(U^n)$ (respectively, $f\in {\mathcal
B}^p_0(U^n)$),
\begin{equation}\label{bou}
\|C_{\phi}f\|_{{\mathcal B}^q}\leq C\|f\|_{\Bl}.
\end{equation}

By Lemma \ref{test} and Inequality (\ref{bou}), there is a $Q\geq0$
such that
$$\sum\limits^n_{k=1}\left|\sum\limits^n_{l=1}\frac{\partial f_w^{(l)}}{\partial
\zeta_l}[\phi(z)] \frac{\partial \phi_l}{\partial
z_k}(z)\right|(1-|z_k|^2)^q\leq CQ$$ for every $w\in U$ and $z\in
U^n.$ Putting $w=\phi_l(z)$ above and using Equations (\ref{par})
and (\ref{par1}), we rewrite the above inequality in the form
$$\sum\limits^n_{k=1}\left|\frac{\partial
\phi_{l}} {\partial z_k}(z)\right|
\frac{\left(1-|z_k|^2\right)^q}{\left(1-|\phi_l(z)|^2\right)^p} \leq
CQ.$$ Condition {\bf (\ref{firsttheorem}c)} follows, and the proof
of Theorem \ref{firsttheorem} is complete. \msk

The paragraph that precedes Corollary \ref{lipresult} justifies
it.
 Therefore, we now proceed to the proof of our next main result:

\medskip

{\it Proof of Theorem \ref{secondtheorem}, Part {\bf (a)}}: By
Theorem \ref{firsttheorem}, $C_\phi$ is a bounded operator from
${\mathcal B}^p(U^n)$ (respectively, ${\mathcal B}^p_0(U^n)$) to
${\mathcal B}^q(U^n)$. By Lemma \ref{phicomp}, $\phi_k\in {\mathcal
B}^q(U^n)$ for each $k\in\{1,2,\ldots,n\}$. The desired statement
then follows from Proposition \ref{precurs506}.

\msk

{\it Proof of Theorem \ref{secondtheorem}, Part {\bf (b)}}: Assume
that $C_{\phi}$ is a compact operator from ${\mathcal B}^p(U^n)$
(respectively, ${\mathcal B}^p_0(U^n)$) to ${\mathcal B}^q(U^n)$, so
that all Conditions {\bf (\ref{firsttheorem}a)} (respectively, {\bf
(\ref{firsttheorem}b)}) and {\bf (\ref{firsttheorem}c)} hold.
Suppose to the contrary that there is a sequence
$(z^{[j]})_{j\in\NN}$ in $U^n$ such that
$\om^{[j]}:=\phi(z^{[j]})\to \pt U^n$ as $j\to\infty$ and that there
is an $\ve>0$ such that for all $j\in\NN$,
\[\sum\limits^n_{k,l=1}\left|\frac{\partial \phi_{l}}{\partial
z_k}(z^{[j]})
\right|\frac{(1-|z^{[j]}_k|^2)^q}{(1-|\phi_l(z^{[j]})|^2)^p} \geq
\ve. \] Surely then, for some $\varepsilon_0>0$ and
$l_0\in\{1,2,\ldots,n\}$, we have that
\[
\sum\limits^n_{k=1}\left|\frac{\partial \phi_{l_0}}{\partial
z_k}(z^{[j]})
\right|\frac{(1-|z^{[j]}_k|^2)^q}{(1-|\phi_{l_0}(z^{[j]})|^2)^p}
\geq \varepsilon_0 \,\,\,\,{\text{ for all }}j\in\NN.
\]
Assume with no loss of generality that $l_0=1$ for notational
convenience; that is, suppose that there is an $\varepsilon_0>0$
such that for all $j\in\NN$,
\begin{equation}\label{hardassump}
\sum\limits^n_{k=1}\left|\frac{\partial \phi_{1}}{\partial
z_k}(z^{[j]})
\right|\frac{(1-|z^{[j]}_k|^2)^q}{(1-|\phi_{1}(z^{[j]})|^2)^p} \geq
\varepsilon_0.
\end{equation}
Since $\phi_1:U^n\rightarrow U$, we must have that
$(\om^{[j]}_1)_{j\in\NN}$ has a subsequence
$(\om^{[j']}_1)_{j\in\NN}$ converging to some $\rho_1\in {\overline
U}$. We consider the two cases $|\rho_1|=1$ and $|\rho_1|<1$; our
approach to both of them is to obtain a contradiction of the
assumption that $C_\phi$ is compact.

{\it Case 1.}  First, suppose that $|\rho_1|=1$.  Then by Lemmas
\ref{precurs3} and \ref{Macla}, $C_\phi$ is not compact.

{\it Case 2.} Second, suppose that $|\rho_1|<1$.  Since
$\om^{[j']}\to\pt U^n$ as $j\rightarrow\infty$, there is an $l_1
\in\{2,3,\ldots,n\}$ such that $|\om^{[j']}_{l_1}|\to 1$ as
$j\to\infty.$ There are two subcases to consider.  Either
\begin{equation}\label{148am}
\lim_{j\to\infty}\left|\frac{\partial \phi_{l_1}}{\partial
z_{k}}(z^{[j']})
\right|\frac{(1-|z^{[j']}_{k}|^2)^q}{(1-|\phi_{l_1}(z^{[j']})|^2)^p}=0\,\,\,\,{\text{
for all }}k\in\{1,2,\ldots,n\},
\end{equation}
or, since Condition ({\bf \ref{firsttheorem}c}) holds, there exists
$k_1\in\{1,2,\ldots,n\}$, $\varepsilon_1>0$, and a further
subsequence $(z^{[j'']})_{j\in\NN}$ of $(z^{[j']})_{j\in\NN}$ such
that
\begin{equation}\label{945pm}
\lim_{j\to\infty}\left|\frac{\partial \phi_{l_1}}{\partial
z_{k_1}}(z^{[j'']})
\right|\frac{(1-|z^{[j'']}_{k_1}|^2)^q}{(1-|\phi_{l_1}(z^{[j'']})|^2)^p}=\varepsilon_1.
\end{equation}

{\it Subcase 2.i:} First, assume the latter subcase above.  Then
there is an $N\in\NN$, such that for all $j\in\NN$ with $j\geq N$,
\[
\left|\frac{\partial \phi_{l_1}}{\partial z_{k_1}}(z^{[j'']})
\right|\frac{(1-|z^{[j'']}_{k_1}|^2)^q}{(1-|\phi_{l_1}(z^{[j'']})|^2)^p}\geq
\frac{\varepsilon_1}{2}.
\]
The quantity on the left side above is bounded independent of $j$,
since $C_\phi$ is bounded by assumption.  It follows from this
fact and the above inequality that for all $j\in\NN$ with $j\geq
N$,
\[
\sum_{k,l=1}^n \left|\frac{\partial \phi_{l_1}}{\partial
z_{k_1}}(z^{[j'']})
\right|\frac{(1-|z^{[j'']}_{k_1}|^2)^q}{(1-|\phi_{l_1}(z^{[j'']})|^2)^p}\geq
\frac{\varepsilon_1}{2}.
\]
Since $(\om^{[j'']}_{l_1})_{j\in\NN}$ is a subsequence of
$(\om^{[j']}_{l_1})_{j\in\NN}$, whose moduli converge to $1$ by
hypothesis in this particular case, we must also have that
$|\om^{[j'']}_{l_1}|\rightarrow 1$ as $j\rightarrow\infty$.
Therefore, the hypotheses of Lemma \ref{precurs3} are satisfied with
``$j$" replaced by ``$\,j''\,$" and ``$\,\varepsilon\,$" replaced by
``$\,\frac{\varepsilon_1}{2}\,$".  By a subsequent application of
Lemma \ref{precurs3} and Lemma \ref{Macla}, we obtain that $C_\phi$
must not be compact.  Thus we have once again obtained a
contradiction.

{\em Subcase 2.ii:} Finally, suppose that Condition (\ref{148am})
holds. We will define a sequence $(h_j)_{j\in\NN}$ in ${\mathcal
B}^p_0(U^n)$ such that is bounded, $h_j\rightarrow 0$ uniformly on
compacta in $U^n$, and $||C_\phi(h_j)||_{{\mathcal B}^p}\not\to0$.
Once again, the compactness of $C_\phi$ will then be contradicted by
Lemma \ref{Macla}, and the proof will be finished.

Define $h_j:\ol{U^n}\rightarrow\CC$ by
\[h_j(z)=h_{\om_{l_1}^{[j']}}^{(l_1)}(z).\]
By Lemma \ref{test}, $(h_j)_{j\in\NN}$ is bounded in ${\mathcal
B}^p(U^{n})$ and also lies in ${\mathcal B}^p_0(U^{n})$.

Let $j\in\NN$, and assume that $K\subset U^n$ is compact. Then there
exists an $r\in (0,1)$ such that $|z_i|\leq r$ for all $z\in K$ and
$i\in\{1,2,\ldots,n\}$.  Hence, for all $z\in K$ and $j\in\NN$,
$$h_j(z)\leq (2+r)\left(\frac{1-|\om^{[j']}_{l_1}|^2}{1-r}\right)^p.$$
It follows that $h_j\to 0$ uniformly on compacta in $U^n$.

We now show that $\left\|C_{\phi}h_j\right\|_{{\mathcal
B}^q}\not\to 0$ as $j\to\infty.$ Now Equations (\ref{hone}) and
(\ref{helle}) imply that
\begin{equation}\label{home1}
\frac{\partial h_{j}}{\partial \zeta_1}\big(\phi(z^{[j']})\big)=
\left(\frac{1-|\om^{[j']}_{l_1}|^2}{1-\om^{[j']}_{l_1}{\overline
\om^{[j']}_{l_1}}}\right)^p =1,
\end{equation}
and
\begin{equation}\label{home2}
\frac{\partial h_{j}}{\partial
\zeta_{l_1}}\big(\phi(z^{[j']})\big)=
\frac{p(2+\om^{[j']}_1)\overline{\om^{[j']}_{l_1}}}{1-|\om^{[j']}_{l_1}|^2}.
\end{equation}
Since $\lim_{j\rightarrow\infty}\om^{[j']}_1=\rho_1$, which is in
$U$ by assumption, there must be an $N'\in\NN$ such that for all
$j\in\NN$ with $j\geq N'$, $|\om_1^{[j']}-\rho_1|<(1-\rho_1)/2$. If
we let $C_{\rho_1}$ be the positive number $1-(1+\rho_1)^2/4$, then
it follows that
\begin{equation}\label{cp1}
1-|\phi_1(z^{[j']})|^2\geq C_{\rho_1}\,\,\,\,{\text{ for all
}}j\in\NN\,\,\,{\text{ such that }}j\geq N'.
\end{equation}
The assumed Condition (\ref{148am}) implies that there is an
$N''\in\NN$ such that for all $j\in\NN$ with $j\geq N''$,
\begin{equation}\label{estimate6}
\sum\limits^n_{k=1}\left|\frac{\partial \phi_{l_1}}{\partial
z_k}(z^{[j']})
\right|\frac{(1-|z^{[j']}_k|^2)^q}{(1-|\phi_{l_1}(z^{[j']})|^2)^p}
\leq \frac{C_{\rho_1}^p\varepsilon_0}{6p}.
\end{equation}
Choose any $j\in\NN$ such that $j\geq\max(N',N'')$, and let
$$I_{j}=\sum\limits^n_{k=1}\left(1-|z^{[j']}_k|^2\right)^q\left|\frac{\partial (h_j^{(l)}\circ\phi)}
{\partial z_k}(z^{[j']})\right|.$$\\
Equation (\ref{h0}) implies that $I_j$ can be written in the form
\[
\sum\limits^n_{k=1}\left(1-|z^{[j']}_k|^2\right)^q\left|\frac{\partial
h_j^{(l)}}{\partial \zeta_1}(\phi(z^{[j']}))\frac{\partial
\phi_1}{\partial z_k}(z^{[j']})+\frac{\partial h_j^{(l)}}{\partial
\zeta_{l_1}}(\phi(z^{[j']}))\frac{\partial \phi_{l_1}}{\partial
z_k}(z^{[j']})\right|.
\]
By Equations (\ref{home1}) and (\ref{home2}), the form of $I_j$
above can be rewritten as
\[
\sum\limits^n_{k=1}(1-|z^{[j']}_k|^2)^q\left|p\frac{(2+\om^{[j']}_1)\overline{\om_{l_1}^{[j']}}}{(1-|\om^{[j']}_{l_1}|^2)}
\frac{\partial \phi_{l_1}}{\partial z_k}(z^{[j']})+\frac{\partial
\phi_1}{\partial z_k}(z^{[j']})\right|,
\]
which is certainly greater than or equal to
\[
\sum\limits^n_{k=1}(1-|z^{[j']}_k|^2)^q\left|\frac{\partial
\phi_1}{\partial
z_k}(z^{[j']})\right|-\sum\limits^n_{k=1}(1-|z^{[j']}_k|^2)^q\left|p\frac{(2+\om^{[j']}_1)
\overline{\om_{l_1}^{[j']}}}{1-|\om^{[j']}_{l_1}|^2} \frac{\partial
\phi_{l_1}}{\partial z_k}(z^{[j']})\right|.
\]
By an elementary estimate and the assumption that $p\geq 1$, it
follows that
\[
I_j\geq\sum\limits^n_{k=1}(1-|z^{[j']}_k|^2)^q\left|\frac{\partial
\phi_1}{\partial z_k}(z^{[j']})\right|
-3p\sum\limits^n_{k=1}\frac{(1-|z^{[j']}_k|^2)^q}{(1-|\om^{[j']}_{l_1}|^2)^p}\left|\frac{\partial
\phi_{l_1}}{\partial z_k}(z^{[j']})\right|
\]
which, by Estimate (\ref{cp1}), is greater than or equal to
\[
C_{\rho_1}^p\sum_{k=1}^n\left|\frac{\partial \phi_1}{\partial
z_{k}}(z^{[j']})
\right|\frac{(1-|z^{[j']}_{k}|^2)^q}{(1-|\phi_1(z^{[j']})|^2)^p}
-3p\sum\limits^n_{k=1}\frac{(1-|z^{[j']}_k|^2)^q}{(1-|\om^{[j']}_{l_1}|^2)^p}\left|\frac{\partial
\phi_{l_1}}{\partial z_k}(z^{[j']})\right|.
\]
By Condition (\ref{hardassump}) and Estimate (\ref{estimate6}), it
follows that for all $j\in\NN$ with $j\geq\max(N',N'')$,
\[
I_j\geq
C_{\rho_1}^p\varepsilon_0-3p\frac{C_{\rho_1}^p\varepsilon_0}{6p}=
\frac{C_{\rho_1}^p\varepsilon_0}{2}>0.
\]
Hence, $||C_\phi h_j||_{\mathcal B}^q\not\to 0$.  The proof of
Theorem \ref{secondtheorem} is now complete.

\medskip

{\it Proof of Corollary \ref{Jac}:} By Equation (\ref{berg}),
Condition {\bf (\ref{Jac}a)} implies that there is a $C\geq 0$ such
that for all $u\in\CC^n$ and $z\in U^n$,
$$\sum_{k=1}^n\frac{|u_k|^2}{(1-|z_k|^2)^2}
\leq C\sum_{l=1}^n\frac{\left|\sum_{k=1}^n\frac{\pt \phi_l}{\pt
z_k}u_k\right|^2}{(1-|\phi_l(z)|^2)^2}.$$ Choosing $u_k=1-|z_k|^2$
for each $k\in\{1,2,\ldots,n\}$, we then obtain that
\begr n&=&\sum_{k=1}^n\frac{(1-|z_k|^2)^2}{(1-|z_k|^2)^2}=\sum_{k=1}^n\frac{|u_k|^2}{(1-|z_k|^2)^2}\nonumber\\
&\leq&C\sum_{l=1}^n\frac{\left(\sum_{k=1}^n\left|\frac{\pt \phi_l}{\pt z_k}\right| |u_k|\right)^2}{(1-|\phi_l(z)|^2)^2}\leq C\left(\sum_{l=1}^n\frac{\sum_{k=1}^n\left|\frac{\pt \phi_l}{\pt z_k}\right| |u_k|}{1-|\phi_l(z)|^2}\right)^2\nonumber\\
&=&C\left(\sum_{k,l=1}^n\left|\frac{\pt \phi_l}{\pt z_k}\right|
\frac{1-|z_k|^2}{1-|\phi_l(z)|^2}\right)^2.\nonumber
\endr
Since $p\geq 1$ and $q\in (0,1]$, we then obtain that
$$\sum_{k,l=1}^n\left|\frac{\pt \phi_l}{\pt z_k}\right| \frac{(1-|z_k|^2)^q}{(1-|\phi_l(z)|^2)^p}\geq\sum_{k,l=1}^n\left|\frac{\pt \phi_l}{\pt z_k}\right| \frac{1-|z_k|^2}{1-|\phi_l(z)|^2}\geq \sqrt {\frac nC},$$
for each $z\in U^n.$  Condition {\bf (\ref{Jac}a)} then follows
from Theorem \ref{secondtheorem}, and the proof of Corollary
\ref{Jac} is complete. \msk

{\em Proof of Corollary \ref{firstiff}:} The equivalence of
Conditions ({\bf \ref{firstiff}a}) and ({\bf \ref{firstiff}b})
immediately follows from Theorem \ref{secondtheorem}.  The fact that
Conditions ({\bf \ref{firstiff}b}) implies Condition ({\bf
\ref{firstiff}c}) follows from Lemma \ref{phicomp}.  Condition {\bf
(\ref{firstiff}c)} implies Condition {\bf (\ref{firstiff}a)} by
Proposition \ref{precurs506}.

\medskip
{\it Proof of Theorem \ref{p unit interval}:} Suppose that
Condition {\bf (\ref{p unit interval}a)} or {\bf (\ref{p unit
interval}b)} holds.  Then $C_{\phi}$ is bounded, and Condition
({\bf \ref{firsttheorem}c}) consequently holds.

We show that Condition ({\bf \ref{p unit interval}c}) holds.  Let
$l\in\{1,2,\ldots,n\}$, and suppose that $(z^{[j]})_{j\in\NN}$ is a
sequence in $U^n$ such that the
 sequence $(\om^{[j]})_{j\in\NN}$ that is given by
$\om^{[j]}=\phi(z^{[j]})$ satisfies
$\lim_{j\to\infty}|\om^{[j]}_l|\to 1.$  For each $j\in\NN$, let
$g_j:{\overline U^n}\rightarrow\CC$ be defined by
\[
g_{j}(z)=g^{(l)}_{\om^{[j]}_l}(z). \] By Lemma \ref{test},
$(g_j)_{j\in\NN}$ is bounded in ${\mathcal B}^p(U^{n})$ and also
lies in ${\mathcal B}^p_0(U^n)$.  An elementary estimate shows
that $g_j\rightarrow 0$ uniformly on compacta in $U^n.$ Therefore,
by Lemma \ref{Macla}, we must have that
$\lim_{j\rightarrow\infty}\left\|C_{\phi}g_{j}\right\|_{{\mathcal
B}^q}=0$.  Hence, for all $j\in\NN$, Equations (\ref{gzero}) and
the chain rule imply that
\begin{eqnarray*}
0&\leq&
\sum\limits^n_{k=1}\frac{(1-|z^{[j]}_k|^2)^q}{(1-|\om^{[j]}_l|^2)^p}
\left|\frac{\partial \phi_l}{\partial z_k}(z^{[j]})\right|\\
&=& \frac{1}{p|\om^{[j]}_l|}\sum\limits^n_{k=1}(1-|z^{[j]}_k|^2)^q
\left|p{\overline
\om^{[j]}_l}\frac{1-|\om^{[j]}_l|^2}{(1-|\om^{[j]}_l|^2)^{p+1}}
\frac{\partial \phi_l}{\partial z_k}(z^{[j]})\right|\\
&=&\frac{1}{p|\om^{[j]}_l|}\sum\limits^n_{k=1}(1-|z^{[j]}_k|^2)^q\left|\frac{\partial
g_j}{\partial \zeta_l}{(\phi(z^{[j]}))}\frac{\partial
\phi_l}{\partial
z_k}(z^{[j]})\right|\\
&=&\frac{1}{p|\om^{[j]}_l|}\sum\limits^n_{k=1}(1-|z^{[j]}_k|^2)^q
\left|\sum_{m=1}^n\frac{\partial g_j}{\partial
\zeta_m}{(\phi(z^{[j]}))}\frac{\partial \phi_m}{\partial
z_k}(z^{[j]})\right|\\
&=&\frac{1}{p|\om^{[j]}_l|}\sum\limits^n_{k=1}(1-|z^{[j]}_k|^2)^q
\left|\frac{\partial(C_\phi(g_j))}{\partial z_k}(z^{[j]})
\right|\\
&\leq&
\frac{1}{p|\om^{[j]}_l|}\left\|C_{\phi}g_{j}\right\|_{{\mathcal
B}^q}\to 0\quad\mbox{as}\quad j\to\infty.
\end{eqnarray*}
Since $l\in\{1,2,\ldots,n\}$ was chosen arbitrarily, we have that
Condition ({\bf \ref{p unit interval}c}) indeed holds.

Next, assume that Condition ({\bf \ref{p unit interval}c}) holds.
Let $(f_j)_{j\in\NN}$ be a sequence in ${\mathcal B}^p(U^n)$
(respectively, ${\mathcal B}^p_0(U^n)$) such that there is a
$C\geq 0$ with $\left\|f_{j}\right\|_{\Bl}\leq C$ for all
$j\in\NN$ and such that $f_{j}\to 0$ uniformly on compacta in
$U^n$. By Lemma \ref{Macla}, it suffices to show that
\begin{equation}\label{939pm}
\lim_{j\to \infty}\left\|C_{\phi}f_{j}\right\|_{{\mathcal B}^q}= 0.
\end{equation}
Once again, if $C=0$, then $f_j=0$ for all $j\in\NN$.  We then have
that $C_\phi(f_j)=0$ for all $j\in\NN$, so that Equation
(\ref{939pm}) holds trivially.  Therefore, we can now assume that
$C>0$.

To conserve space, for $k,l\in\{1,2\ldots,n\}$ we define
$B^\phi_{k,l}:U^n\rightarrow\RR$ by
\[
B^\phi_{k,l}(z)=\left|\frac{\partial \phi_l} {\partial
z_k}(z)\right| (1-|z_k|^2)^q.
\]
Now let $\varepsilon>0$.  Since $f_j\rightarrow 0$ on compacta and
in particular on the compact set $\{\phi(0)\}$, there is an
$N\in\NN$ such that
\begin{equation}\label{1145}
|f_j[\phi(0)]|<\frac{\varepsilon}{2}\,\,\,\,{\text{ for all
}}\,\,\,j\in \NN\,\,{\text{ such that}}\,\,j\geq N.
\end{equation}
Condition ({\bf \ref{p unit interval}c}) implies that there exists
an $r\in(0,1)$ such that whenever $l\in\{1,2,\ldots,n\}$ and
$|\phi_l(z)|>r,$
\begin{equation}\label{930pm}
\sum\limits^n_{k=1}\left|\frac{\partial \phi_l}{\partial
z_k}(z)\right|\frac{(1-|z_k|^2)^q}{(1-|\phi_l(z)|^2)^p}<\frac{\varepsilon}{2n^2C}.
\end{equation}
For $l\in\{1,2,\ldots,n\}$, let
\[
F^{(l)}_r=\{z\in U^n\, |\, |\phi_l(z)|> r\} \] and
\[G_r^{(l)}=U^n-F^{(l)}_r=\{z\in U^n:|\phi_l(z)|\leq r\}.
\]
If $k,l\in\{1,2,\ldots,n\}$, $j\in\NN$, $z\in U^n$, and $z\in
F^{(l)}_r$, then
\begin{eqnarray}
\left|\frac{\partial f_j} {\partial w_l}(\phi(z))\right|
\left|\frac{\partial \phi_l} {\partial z_k}(z)\right|
(1-|z_k|^2)^q &\leq & \|f_j\|_p\left|\frac{\partial
\phi_l}{\partial
z_k}(z)\right|\frac{(1-|z_k|^2)^q}{(1-|\phi_l(z)|^2)^p}\nonumber\\
&<&C\frac{\varepsilon}{2n^2C}\nonumber\\
&=&\frac{\varepsilon}{2n^2}\label{didnt}.
\end{eqnarray}
We now obtain the same estimate for $z\in G^{(l)}_r$.  By Lemma
\ref{phicomp}, we have that $\phi_l\in {\mathcal B}^q(U^n)$ for all
$l\in\{1,2,\ldots,n\}$.  It follows that there is an $M\geq 0$ such
that for all $k,l\in\{1,2,\ldots,n\}$ and $z\in U^n$,
\begin{equation}\label{bphi}
B^\phi_{k,l}(z)\leq M.
\end{equation}
Let $z\in U^n$, $k, l\in\{1,2,\ldots,n\}$, and $j\in\NN$, and
define $g_{z,j}^{(l)}:U\rightarrow \CC$ by
$g_{z,j}^{(l)}(\xi)=f_j(\xi,z_l')$. By Lemma \ref{hln1},
$f_j\in{\mathcal L}^{1-p}(U^n)$, and $f_j$ extends continuously to
${\overline U^n}$ by Lemma \ref{Cauchy}. If we once more label the
continuous extension of $f_j$ by $f_j$, then by Lemma
\ref{usedonce?}, there is an $N'$ such that
\[|f_j(z)|<\frac{r\varepsilon}{2n^2M}\,\,{\text{ for all
}}\,z\in{\overline U^n}\,{\text{ and }}\,j\in\NN\,{\text{ such that
}}\,j\geq N'.
\]
Furthermore, $g_{z,j}^{(l)}$ is bounded on $B_r:=\{\xi\in
U:|\xi|\leq r\}$.  From these two facts and Cauchy's estimate, we
have that for all $\xi\in B_r$,
\[\left|\left(\frac{d}{d\xi}g_{z,j}^{(l)}\right)(\xi)\right|\leq
\frac{\sup_{\nu\in
B_{(1+r)/2}}\left|g_{z,j}^{(l)}(\nu)\right|}{\left(\frac{1+r}{2}\right)}=
\frac{2}{1+r}\sup_{\nu\in
B_{(1+r)/2}}\left|g_{z,j}^{(l)}(\nu)\right|.
\]
Let $j\in\NN$, $j\geq N'$, and $z\in G_r^{(l)}$. Then
$\phi_l(z)\in B_r$, and we have that
\[
\left|\frac{\partial f_j} {\partial \zeta_l}(\phi(z))\right|
=\left|\left(\frac{d}{d\xi}g_{z,j}^{(l)}\right)(\phi_l(z))\right|
\leq \frac{2}{1+r}\sup_{\nu\in
B_{(1+r)/2}}\left|g_{z,j}^{(l)}(\nu)\right|
\]
\begin{eqnarray*}
&=&\frac{2}{1+r}\sup_{\nu\in
B_{(1+r)/2}}\left|f_j(\nu,z_l')\right|\\
&\leq & \frac{2}{1+r}\sup_{z\in{\overline U^n}}\left|f_j(z)\right|\\
&<&
\frac{2}{1+r}\left[\frac{r\varepsilon}{2n^2M}\right]<\frac{\varepsilon}{2n^2M}.
\end{eqnarray*}
It follows from Inequality (\ref{bphi}) that for all $z\in
G^{(l)}_r$ and $j\in\NN$ such that $j\geq N'$,
\[
\left|\frac{\partial f_j} {\partial
\xi_l}[\phi(z)]\right|B^\phi_{k,l}(z)<\frac{\varepsilon}{2n^2M}\cdot
 M=\frac{\varepsilon}{2n^2}.
\]
Now $U^n=F^{(l)}_r\cup G^{(l)}_r$, and Estimate (\ref{didnt})
implies that
\begin{equation}\label{estim}
\left|\frac{\partial f_j} {\partial
\xi_l}[\phi(z)]\right|B^\phi_{k,l}(z)<\frac{\varepsilon}{2n^2}
\end{equation}
for $z\in F^{(l)}_r$ and $j\in\NN$.  It follows that Inequality
(\ref{estim}) holds for all $z\in U^n$, $k,l\in\{1,2\ldots,n\}$, and
$j\in\NN$ such that $j\geq N'$.  Letting $N''=\max(N,N')$, we
conclude from the above estimate and Estimate (\ref{1145}) that for
all $z\in\NN$ and $j\in\NN$ with $j\geq N''$,
\[
|[C_\phi(f_j)](0)|+\sum_{k=1}^n\left|\frac{\partial (C_{\phi}f_j)}
{\partial z_k}(z)\right| (1-|z_k|^2)^q
\]
\begin{eqnarray*}
 &\leq
&\frac{\varepsilon}{2}+\sum_{k=1}^n\sum\limits^n_{l=1}
\left|\frac{\partial f_j} {\partial
w_l}(\phi(z))\right| B^\phi_{k,l}(z)\\
&\leq &
\frac{\varepsilon}{2}+\sum_{k=1}^n\sum\limits^n_{l=1}\frac{\varepsilon}{2n^2}\\
&=&\varepsilon.
\end{eqnarray*}
Therefore, Equation (\ref{939pm}) holds, and Conditions {\bf
(\ref{p unit interval}a)} and {\bf (\ref{p unit interval}b)}
follow.  The proof of Theorem \ref{p unit interval} is now
complete.

\medskip

{\it Proof of Corollary \ref{cor of t2and3}:} Let $q>0$.  By
assumption, if $p\geq 1$ (resp., $p\in(0,1)$), then there is an
$M\geq 0$ such that for all $l\in\{1,2,\ldots,n\}$ and $z\in U^n$,
$$\sum\limits^n_{k,l=1}\left|\frac{\partial \phi_{l}}
{\partial
z_k}(z)\right|\frac{(1-|z_k|^2)^q}{(1-|\phi_l(z)|^2)^p}\leq
\sum\limits^n_{l=1}\frac{\|\phi_l\|_q}{(1-\|\phi_l\|_\infty^2)^p}\leq
M.$$ Hence by Theorem \ref{firsttheorem}, $C_{\phi}$ is a bounded
operator from ${\mathcal B}^p(U^n)$ and ${\mathcal B}^p_0(U^n)$ to
${\mathcal B}^q(U^n)$.
 The assumption that $\|\phi_l\|_\infty<1$ for each $l\in\{1,2,\ldots,n\}$ implies that
 Condition (\ref{firstcom})
(respectively, Condition ({\bf \ref{p unit interval}c})) is
vacuously satisfied if $p\geq 1$ (respectively, $p\in (0,1)$).
The desired statement then follows from a simultaneous application
of Theorems \ref{secondtheorem} and Theorem \ref{p unit interval}.

\medskip

{\it Proof of Corollary \ref{Stevo}:} By Lemma
\ref{Schwarzlemmapolydisk}, there is a positive constant $C>0$ such
that
$$H_{\phi(z)}\big(J\phi(z)u,\ol{J\phi(z)u}\big)\leq CH_z(u,\bar u)$$
for every $z\in U^n$ and $u\in \CC^n$.  By (\ref{berg}), we then
have that
$$\sum_{l=1}^n\frac{\left|\sum_{k=1}^n\frac{\pt \phi_l}{\pt z_k}u_k\right|^2}{(1-|\phi_l(z)|^2)^2}
\leq C\sum_{k=1}^n\frac{|u_k|^2}{(1-|z_k|^2)^2}.$$ Successively
$u$ in the above inequality by
$$u^{(k)}:=(\underbrace{0,\ldots,0}_{k-1},1-|z_k|^2,\underbrace{0,\ldots,0}_{n-k})=([1-|z_k|^2],0'_k)$$
for each $k\in\{1,2,\ldots,n\}$, we obtain the inequality
$$\sum_{l=1}^n\frac{\left|\frac{\pt \phi_l}{\pt z_k}(1-|z_k|^2)\right|^2}{(1-|\phi_l(z)|^2)^2}
\leq C$$
for each such $k$ and all $z\in U^n$.  It follows that
$$\left|\frac{\pt \phi_l}{\pt z_k}(z)\right|\frac{1-|z_k|^2}{1-|\phi_l(z)|^2}\leq C$$
for each $k,l\in \{1,\ldots,n\}$ and $z\in U^n$.

Using the above inequality and the assumptions that $q\geq 1$ and
$p\in (0,1)$, we obtain that for all $z\in U^n$,
\begin{eqnarray}
\sum\limits^n_{k=1}\left|\frac{\partial \phi_{l}} {\partial
z_k}(z)\right|\frac{(1-|z_k|^2)^q}{(1-|\phi_l(z)|^2)^p}
&\leq& C\sum\limits^n_{k=1}(1-|z_k|^2)^{q-1}(1-|\phi_l(z)|^2)^{1-p}\nonumber\\
&\leq&C\sum\limits^n_{k=1}(1-|\phi_l(z)|^2)^{1-p}\label{stran}\\
&\leq& Cn.\nonumber
\end{eqnarray}

It follows from Theorem \ref{firsttheorem} that $C_{\phi}$ is a
bounded operator from ${\mathcal B}^p(U^n)$ and ${\mathcal
B}^p_0(U^n)$ to ${\mathcal B}^q(U^n).$  Letting $|\phi_l(z)|\to 1$
in Inequality (\ref{stran}), we see that Condition ({\bf \ref{p
unit interval}c}) holds. Therefore, by Theorem \ref{p unit
interval}, $C_{\phi}$ is a compact operator from ${\mathcal
B}^p(U^n)$ and ${\mathcal B}^p_0(U^n)$ to ${\mathcal B}^q(U^n)$.
The proof is now complete.

\section{Composition operators between ${\mathcal B}^p_0(U^n)$ and ${\mathcal B}^q_0(U^n)$}\label{little spaces}

We now formulate and prove the corresponding results concerning
$C_\phi$ when it maps between little $p$- and $q$-Bloch spaces of
$U^n$.  First, we prove the following boundedness result
corresponding to Theorem \ref{firsttheorem}:

\begin{theorem}\label{almostmissed} Let $\phi$ be a holomorphic self-map of $U^n$ and assume that $p,q>0$.
Then $C_\phi$ is a bounded operator from ${\mathcal B}_0^p(U^n)$ to
${\mathcal B}_0^q(U^n)$ if and only if

{\bf (\ref{almostmissed}a)} $\phi^\gamma\in {\mathcal B}_0^q(U^n)$
for every $n$-multi-index $\gamma$, and

{\bf (\ref{almostmissed}b)} there exists a constant $C$ such that
for all $z\in U^n$,
\[
\sum\limits^n_{k,l=1}\left|\frac{\partial \phi_{l}} {\partial
z_k}(z)\right|\frac{\left(1-|z_k|^2\right)^q}
{\left(1-|\phi_l(z)|^2\right)^p}\leq C.
\]
\end{theorem}
\begin{proof}$\Rightarrow$: Suppose that $C_\phi$ is a bounded operator from
${\mathcal B}_0^p(U^n)$ to ${\mathcal B}_0^q(U^n)$. Then Condition
{\bf (\ref{firsttheorem}b)} holds.  Condition {\bf
(\ref{almostmissed}b)} in turn follows from Theorem
\ref{firsttheorem}. Also, $\phi^\gamma\in {\mathcal B}_0^q(U^n)$
by Lemma \ref{phicomp}.

$\Leftarrow$: By Theorem \ref{firsttheorem}, $C_{\phi}$ is a
bounded operator from ${\mathcal B}^p(U^n)$ to ${\mathcal
B}^q(U^n)$. Denoting the norm of $C_\phi$ by $||C_\phi||_{p,q}$,
we then have that
\[
\|C_\phi f\|_{{\mathcal B}^q}\leq
||C_\phi||_{p,q}\,\|f\|_{\Bl}\quad \mbox{for all}\; f\in{\mathcal
B}_0^p(U^n).
\]
Therefore, it remains to show that $C_\phi f\in{\mathcal
B}_0^q(U^n)$ if $f\in {\mathcal B}_0^p(U^n).$ Let $\varepsilon>0$,
and suppose that $f\in {\mathcal B}_0^p(U^n)$. Then there is a
polynomial function $p_\ve:\CC^n\rightarrow\CC$ such that
\begin{equation}\label{tonight}
\|f-p_\ve\|_{\Bl}<\frac{\ve}{2||C_\phi||_{p,q}}.
\end{equation}
For any $n$-multi-index $\gamma$, define
$p_\gamma:\CC^n\rightarrow\CC$ by $p_\gamma(z)=z^\gamma$. The
assumption that $\phi^\gamma\in {\mathcal B}_0^q(U^n)$ and the
identity $\phi^\gamma=C_\phi(p_\gamma)$ for every $n$-multi-index
$\gamma$ together imply that $C_\phi(p_\ve)\in{{\mathcal
B}_0^q}(U^n)$.  In turn, there is a polynomial
$p'_\ve:\CC^n\rightarrow\CC$ such that
$||C_\phi(p_\ve)-p'_\ve||_{{\mathcal B}^q}<\ve/2$. Therefore, by
Inequality (\ref{tonight}), we have that
\begin{eqnarray*}
||C_\phi(f)-p'_\ve||_{{\mathcal B}^q} &\leq
&||C_\phi(f)-C_\phi(p_\ve)||_{{\mathcal
B}_q}+||C_\phi(p_\ve)-p'_\ve||_{{\mathcal B}_q}\\
&\leq&||C_\phi||_p^q||f-p_\ve||_{{\mathcal
B}_q}+||C_\phi(p_\ve)-p'_\ve||_{{\mathcal B}^q}\\
&<&\frac{\ve}{2}+\frac{\ve}{2}=\varepsilon.
\end{eqnarray*}
Therefore, $C_\phi(f)\in{\mathcal B}_0^q(U^n)$.
\end{proof}
Prior to our list of compactness results between ${\mathcal
B}^p_0(U^n)$ and ${\mathcal B}^q_0(U^n)$, we state and prove the
following prerequisite fact:
\begin{lemma}\label{littlelemma}
Suppose that $p,q>0$, and assume that $\phi:U^n\rightarrow U^n$ is
holomorphic.  Let $C_\phi$ be a compact operator from ${\mathcal
B}^p_0(U^n)$ to ${\mathcal B}^q(U^n)$, and suppose that
$\phi^\g\in{\mathcal B}^q_0(U^n)$ for every $n$-multi-index $\g$.
Then $C_\phi$ is more particulary a compact operator from ${\mathcal
B}^p_0(U^n)$ to ${\mathcal B}^q_0(U^n)$.
\end{lemma}
\begin{proof}
Let $f\in{\mathcal B}_0^p(U^n)$ and $r\in (0,1)$.  Define
$f_r:\frac{1}{r}{\overline U^n}\rightarrow\CC$ by
 $f_r(z)=f(rz)$.  If $z\in U^n$, then
\begr 0& \leq &|f(0)|+\sum\limits^n_{k=1} \left|\frac{\partial
f_r}
{\partial z_k}(z)\right|\left(1- |z_k|^2\right)^p\nonumber\\
&=&|f(0)|+\sum\limits^n_{k=1}
r\left|\frac{\partial f} {\partial z_k}(rz)\right|\left(1- |z_k|^2\right)^p\nonumber\\
&\leq&|f(0)|+\sum\limits^n_{k=1} \left|\frac{\partial f} {\partial
z_k}(rz)\right|\left(1-
|rz_k|^2\right)^p\nonumber \\
&\leq&\|f\|_{\Bl}<\infty.\nonumber
\endr
It follows that $f_r\in{\mathcal B}^p(U^n)$. Since $f_r\in H(\frac
1rU^n)$, there is a polynomial $p_r:\CC^n\rightarrow\CC$ such that
for all $z\in \frac{2}{1+r}{\overline U^n}$,
\begin{equation}\label{getgoing}
|f_r(z)-p_r(z)|<(1-r)^2.
\end{equation}
For all $z\in U^n$ and $k\in \{1,2,\ldots,n\}$ we then have that
\begin{eqnarray*}
\left|\frac{\pt f_r}{\pt z_k}(z)-\frac{\pt p_r}{\pt
z_k}(z)\right|&\leq &\frac 1{2\pi}\left|\int_{\frac2{1+r}\pt
U}\frac{f_r(\zeta_k,z_k')-p_r(\zeta_k,z_k')}{(z_k-\zeta_k)^2}d\zeta_k\right|\\
&\leq& 8.
\end{eqnarray*}
Therefore, $f_r-p_r\in {\mathcal B}^p(U^n)$; in fact,
$||f_r-p_r||_{{\mathcal B}^p}\leq 8$.  It follows that
\[
||p_r||_{{\mathcal B}^p}\leq ||f_r||_{{\mathcal B}^p}+8.
\]
Let $r_j$ be a sequence in $(0,1)$ such that
$\lim_{j\to\infty}r_j=1.$  Inequality (\ref{getgoing}) shows that
$p_{r_j}-f\to 0$ uniformly on compacta in $U^n$; furthermore, the
sequence $(p_{r_j}-f)_{j\in\NN}$ is norm-bounded (by $8$) in
${\mathcal B}^p(U^n)$. Since $C_\phi$ is compact from ${\mathcal
B}^p(U^n)$ to ${\mathcal B}^q(U^n)$, Lemma \ref{Macla} implies that
\begin{equation}\label{yahooalmostdone}
\|C_\phi p_{r_j}-C_\phi f\|_{{\mathcal
B}^q}=\|C_\phi(p_{r_j}-f)\|_{{\mathcal B}^q}\to 0\,\,\,{\text{ as
}}\,\,j\to\infty.
\end{equation}
By the assumption that $\phi^\g\in {\mathcal B}_0^q(U^n)$ for every
$n$-multi-index $\g$ and the linearity of $C_\phi$, we have that
$C_\phi p_{r_j}\in {\mathcal B}_0^q(U^n)$ for all $j\in\NN$.
Combining this fact with Equation (\ref{yahooalmostdone}), we obtain
that $C_\phi(f)\in{\mathcal B}^q_0(U^n)$.  Since $f\in{\mathcal
B}^p_0(U^n)$ was chosen arbitarily, the proof of the lemma is now
complete.
\end{proof}

Finally, we state and prove the following results concerning compact
composition operators between little $p$- and little $q$-Bloch
spaces of $U^n$:
\begin{theorem}\label{littletheorem} Let $\phi$ be a holomorphic self-map of
$U^n$, and suppose that $p, q>0.$

\item[{\bf \,\,\,(a)}] If $C_\phi$ is a bounded operator from
${\mathcal B}^p_0(U^n)$ to ${\mathcal B}^p_0(U^n)$ and
\begin{equation}\label{littlecom}
\lim_{\phi(z)\to \partial
U^n}\sum\limits^n_{k,l=1}\left|\frac{\partial \phi_{l}} {\partial
z_k}(z)\right|\frac{(1-|z_k|^2)^q}{(1-|\phi_l(z)|^2)^p} =o(1),
\end{equation}
then $C_{\phi}$ is a compact operator from ${\mathcal B}_0^p(U^n)$
to ${\mathcal B}_0^q(U^n)$.

\item[\,\,\,{\bf (b)}] Let $p\geq 1$ and $q>0.$ If $C_{\phi}$ is a
compact operator from ${\mathcal B}_0^p(U^n)$ to ${\mathcal
B}_0^q(U^n)$, then $\phi^\gamma\in {\mathcal B}_0^q(U^n)$ for
every $n$-multi-index $\gamma$, and Condition $(\ref{littlecom})$
holds.
\end{theorem}

\begin{proof}{\bf (a)}: $C_{\phi}$ is compact
between ${\mathcal B}^p(U^n)$ and ${\mathcal B}^q(U^n)$ by Theorem
\ref{secondtheorem}.  The statement {\bf (a)} then follows from
Lemma \ref{littlelemma}.

{\bf (b)}: Let $\gamma$ be an $n$-multi-index, and define the
polynomial function $p^\g:\CC^n\rightarrow\CC$ by
$p_\gamma(z)=z^\gamma$.  Since $p_\gamma\in{\mathcal B}_0^p(U^n)$,
we obtain by hypothesis that $C_\phi(p_\g)=\phi^\g \in {\mathcal
B}_0^q(U^n)$. Since $C_\phi$ maps ${\mathcal B}^p_0(U^n)$ to
${\mathcal B}^q_0(U^n)$ by assumption, $C_\phi$ certainly maps
${\mathcal B}^p_0(U^n)$ to ${\mathcal B}^q(U^n)$.  Equation
(\ref{littlecom}) then follows from Theorem \ref{secondtheorem}.
\end{proof}
In direct correspondence to Corollary \ref{firstiff}, we now have
the following characterization of holomorphic self-maps of $U^n$
that induce compact composition operators between little $p$- and
little $q$-Bloch spaces of $U^n$ in the case that $p\geq 1$.

\begin{corollary} Let $\phi$ be a holomorphic self-map of $U^n$, and
assume that $p\geq 1$ and $q>0.$  Then $C_{\phi}$ is a compact
operator from $ {\mathcal B}_0^p(U^n)$ to ${\mathcal B}_0^q(U^n)$
if and only if $\phi^\g\in {\mathcal B}_0^q(U^n)$ for every
$n$-multi-index $\g$ and Equation $(\ref{littlecom})$ holds.
\end{corollary}

Corresponding to Corollary \ref{p unit interval} is the following
result, which resolves the issue of compactness of composition
operators between little $p$- and little $q$-Bloch spaces of $U^n$
in the case that $p\in(0,1)$:
\begin{theorem}\label{littletheorem2} Let $\phi$ be a holomorphic self-map
of $U^n$, and suppose that $p\in (0,1)$ and $q>0.$ Then $C_{\phi}$
is a compact operator from ${\mathcal B}_0^p(U^n)$ to ${\mathcal
B}_0^q(U^n)$ if and only if $\phi^\g\in {\mathcal B}_0^q(U^n)$ for
every multi-index $\g$, $C_\phi$ is a bounded operator from
${\mathcal B}_0^p(U^n)$ to ${\mathcal B}_0^q(U^n)$, and for each
$l\in\{1,2,\ldots,n\}$, \bege\label{almostalmost}
\lim_{|\phi_l(z)|\to 1}\sum\limits^n_{k=1}\left|\frac{\partial
\phi_{l}} {\partial
z_k}(z)\right|\frac{\left(1-|z_k|^2\right)^q}{\left(1-|\phi_l(z)|^2\right)^p}=0.
\ende
\end{theorem}

\begin{proof}$\Rightarrow$: Once again, $\phi^\g\in {\mathcal B}_0^q(U^n)$
for every multi-index $\g$ by Lemma \ref{phicomp}. Since $C_{\phi}$
is a compact operator from ${\mathcal B}_0^p(U^n)$ to ${\mathcal
B}_0^q(U^n)$, then it obviously maps into ${\mathcal B}^q(U^n)$ and
certainly compact as an operator from ${\mathcal B}_0^p(U^n)$ to
${\mathcal B}^q(U^n)$.  Therefore, Equation (\ref{almostalmost})
holds by Theorem \ref{p unit interval}.

$\Leftarrow$: Surely, $C_\phi$ is a bounded operator from ${\mathcal
B}_0^p(U^n)$ to ${\mathcal B}^q(U^n)$.  $C_\phi$ is furthermore
compact between these two spaces by Theorem \ref{p unit interval}.
The compactness of $C_\phi$ from ${\mathcal B}^p_0(U^n)$ to
${\mathcal B}^q_0(U^n)$ then follows from Lemma \ref{littlelemma}.
\end{proof}

Finally, we state and prove the following corollary:
\begin{corollary} Suppose that $p,q>0$, and let $\phi$ be a holomorphic self-map of $U^n$ such that
$\|\phi_l\|_\infty<1$ for each $l\in\{1,2,\ldots,n\}$, with
$\phi^\g\in {\mathcal B}_0^q(U^n)$ for every multi-index $\g$.  Then
$C_\phi$ is a compact operator from ${\mathcal B}_0^p(U^n)$ to
${\mathcal B}_0^q(U^n)$.
\end{corollary}
\begin{proof}
If $\|\phi_l\|_\infty<1$ for each $l\in\{1,2,\ldots,n\}$, then
Conditions (\ref{littlecom}) and (\ref{almostalmost}) are vacuously
satisfied.  For $p\geq 1$, the statement of the corollary then
follows from Theorem \ref{littletheorem}, and for $p\in(0,1)$, the
statement follows from Theorem \ref{littletheorem2}.
\end{proof}
\section{An Open Problem}

For the case that $p,q\in(0,1)$, we conjecture that compactness of
$C_\phi$ from ${\mathcal B}^p(U^n)$ or ${\mathcal B}^p_0(U^n)$ to
${\mathcal B}^q(U^n)$ implies that $||\phi_l||_\infty<1$ for all
$l\in\{1,2,\ldots,n\}$. This conjecture is known to hold when
$p=q$ and $U^n$ is replaced by $B_n$
\cite[~Ch.~4]{Cowen/MacCluer}.
\section*{Acknowledgements}
Thanks are extended to Joe Ball, Hari Bercovici and Richard Timoney
for helpful comments on an earlier version of the manuscript.

\end{document}